\documentclass[final,10pt]{elsarticle}

\usepackage{amssymb,amsmath,amsthm,mathtools,dsfont,centernot}

\def\DZ              {\texttt{DZ}}
\newcommand{\HN}{H_N}       
\newcommand{\HIE}{H_{\IE}}  
\newcommand{\ftrue}{f_{\IE}}
\newcommand{\Avg}{\IE_N}

\usepackage{bm}                 
\DeclareMathOperator*{\argmin}{\textrm{argmin}\,} 
\def\Diag            {\textrm{Diag}} 	 
\def\dist            {\textrm{dist}}     
\def\nnz             {\textrm{nnz}}      
\def\ones            {\bm{1}}            
\def\rk              {\textrm{rank}}     
\def\st              {\,\textrm{s.t.}\,} 
\def\zeros           {\bm{0}}            

\def\IE {\mathbb{E}} 		
\def\II {\mathbb{I}} 		%
\def\IP {\mathbb{P}} 		
\def\IR  {\mathbb{R}} 	    
\def\M {\bm{M}}

\def\U {\bm{U}}

\def\cB { \mathcal{B}}

\def\cD { \mathcal{D}}
\def\cE { \mathcal{E}}

\def\cK { \mathcal{K}}

\def\cN { \mathcal{N}}
\def\cO { \mathcal{O}}
\def\cP { \mathcal{P}}

\def\cR { \mathcal{R}}

\def\cU { \mathcal{U}}
\def\cV { \mathcal{V}}

\def\cX { \mathcal{X}}

\def\bPi      {\boldsymbol{\Pi}}

\usepackage{tikz, pgfplots}\pgfplotsset{compat=1.17}
\usepackage{xcolor}

\newtheorem{theorem}{Theorem}[section] 
\newtheorem{remark}[theorem]{Remark} 
\newtheorem{proposition}[theorem]{Proposition}
\newtheorem{corollary}[theorem]{Corollary} 
\newtheorem{lemma}[theorem]{Lemma}
\newtheorem{definition}[theorem]{Definition}

\usepackage{tabularx} 
\usepackage{enumitem} 
\usepackage{multirow,multicol}
\usepackage{makecell}
\usepackage{comment}
\usepackage{url}
\usepackage{booktabs}
\usepackage{tcolorbox} 
\tcbset{left=1pt, right=1pt, top=1pt, bottom=1pt, boxsep=0pt, sharp corners,
 colback=gray!9, 
 colframe=black, 
 boxrule=0pt, 
 before skip=0pt, 
 after skip=0pt,
}

\usepackage{caption}
\newcommand{\yes}{\checkmark}
\newcommand{\no}{} 

\usepackage[ruled,vlined]{algorithm2e}
\SetKwInput{KwInput}{Input} 
\SetKwInput{KwOutput}{Output} 
\DontPrintSemicolon

\usepackage{cleveref}
\crefname{section}{§}{§§}
\Crefname{section}{§}{§§}
 
\usepackage{etoc} 

\journal{Artificial Intelligence}
\begin{document}
\begin{frontmatter}
\title{Optimal Network Pricing for Oblivious Users under Projected Decision-Dependent Distributions
}
\author{\vspace{-4mm} Yixuan Li\,$^\star$\tnoteref{prelimnote}, ~~~Andersen Ang\,$^\star$, ~~~Sebastian Stein\vspace{-1mm}
\tnotetext[prelimnote]{Equal contribution. 
A preliminary arXiv preprint: \url{https://arxiv.org/pdf/2510.07157}.
}
}
\affiliation{organization={School of Electronics and Computer Science, University of Southampton},country={United Kingdom}
}

\begin{abstract}
Efficient large-scale network allocation requires data-driven pricing mechanisms that internalize stochastic, nonlinear user behavior. 
We move beyond the classic fully strategic agents to study oblivious users (agents with bounded rationality and imperfect information).
Rather than assuming an infinite horizon, our regime acknowledges that real-world flows are too transient to equilibrate among users.
We introduce a novel Optimal Network Pricing (ONP) problem for such users, which induces Performativity: a Decision-Dependent environment where pricing decisions endogenously shift the flow distribution.
Without a closed-form distribution, the platform must learn optimal prices from sampled responses.

This setting introduces a new challenge: capacity boundaries and projection operators make the optimization landscape nonsmooth, invalidating gradient-based methods.
We show that a widely adopted optimality concept Performative Stability (PS) fails in ONP, collapsing to a trivial solution.
We then define a new optimality concept, the Projected Performative Optimum ($\Pi\text{PO}$) for the unique global optimum.
Targeting $\Pi\text{PO}$ is algorithmically hard given the performative nonsmooth Jacobian,
so we propose a novel framework combining Sample Average Approximation with Trust-Region Sequential Quadratic Programming, explicitly handling the capacity boundaries,
with theoretical guarantees on probabilistic convexity, sample complexity, and computational complexity.
Experiments show that our $\Pi$PO solver significantly outperforms PS-seeking heuristics and a proposed baseline (improving social welfare by 81\% on GÉANT), highlighting that properly handling capacity boundaries unlocks substantial gains in social welfare. 
More broadly, this work advances intelligent systems that learn under performative, capacity-constrained feedback, a core challenge in real-world AI applications.
\end{abstract}
{
}
\end{frontmatter}

\localtableofcontents

\begin{figure}[h!]
\centering\includegraphics[width=0.99\linewidth]{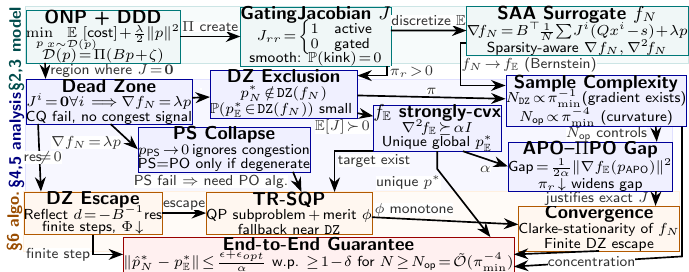}
\caption{The flow of the paper and the relationship.}
\label{fig:organization}
\end{figure}

\section{Introduction}
Efficiently optimizing large-scale network systems remains challenging, due to the tension between finite network capacity\footnote{E.g., road space, bandwidth, computational cores.} and stochastic, often unpredictable demand.
Modern infrastructures increasingly use dynamic pricing to mediate this tension: smart grids deploy time-of-use tariffs for load shifting \cite{siano2014demand}; communication networks use variable pricing to mitigate congestion and optimize spectrum resources \cite{sen2013smart};
Payment Channel Networks \cite{sivaraman2018routing} use variable pricing to ensure quality of service.
Across these domains, pricing serves as a control signal shaping congestion patterns\footnote{Such as smoothing peak demand to reduce congestion.} and allocating resources efficiently \cite[Ch4]{small2024economics}.
Yet setting this signal well, increasingly requires data-driven, learning-based mechanisms that account for how pricing endogenously reshapes the flow distribution, a response that depends on how users behave.

\paragraph{Strategic Agent}
The efficacy of pricing mechanisms depends on behavioral assumptions about the users they seek to influence.
Network flow modeling has historically focused on equilibrium-based paradigms\footnote{Wardrop \cite{wardrop1952road}, Rosenthal's congestion games \cite{rosenthal1973class}. }, relying on assumption A0:
\begin{tcolorbox}
\underline{\textbf{(A0) Strategic Agents}} are hyper-rational, omniscient users with \textit{perfect information and complete understanding} of network topology, instantaneous link congestion, and all other agents' decisions. 
\end{tcolorbox}
\vspace{1mm}
\noindent While theoretically elegant, these assumptions are rarely met\footnote{
E.g., cloud tenants rarely monitor millisecond-level load changes in a provider's internal backbone; drivers lack real-time visibility into citywide traffic conditions.
}, 
which constitutes a critical gap for intelligent pricing systems that must operate under real human behavior.
Real-world users instead exhibit \emph{bounded rationality} \cite{simon1955behavioral}, acting on habit or transient preferences, and \emph{imperfect information}, since platforms seldom disclose other users' real-time behavior, and users' limited attention hinders instantaneous processing of all available information.

A further obstacle is temporal: classical equilibrium-based control assumes an \emph{infinite horizon}, but real systems usually do not stand still long enough to equilibrate among users \cite{blum2006routing, young2004strategic}.
In practice, users are too transient and dynamic to settle into such a equilibrium state. 
Within a finite window, they lack the repeated interactions needed to \emph{learn} true network conditions,
so this incomplete learning emerges as stochastic behavior rather than user equilibrium. 
Moreover, constant user turnover compounds this: as the "set of players" is always changing, the system never reaches a user equilibrium, which requires capturing a snapshot of the transient dynamics.

\paragraph{Oblivious User}
To address A0's shortcomings, we adopt a more empirically grounded behavioral assumption (A1).
\begin{tcolorbox}
\underline{\textbf{(A1) Oblivious users} \cite{karakostas2011degradation}}
choose routes based on fixed criteria (e.g.~posted prices and static distance) or transient preferences (e.g. willingness to join temporary events), oblivious to other users' momentary behavior.
\end{tcolorbox}
\noindent Under A1, human decisions are stochastic, responding to aggregate signals rather than granular state data \cite{karakostas2011degradation}.
It defines a distinct yet practically pervasive regime\footnote{
E.g., in cloud computing, tenants choose availability zones and reserve instances at posted prices based on advertised or historical load, yet cannot observe other tenants' concurrent decisions, and so cannot anticipate the resulting transient spikes; in transportation, commuters follow habitual routes guided by predicted congestion, yet remain blind to other drivers' instantaneous choices and the resulting momentary flow fluctuations.}. 
Crucially, obliviousness does not mean total ignorance: users may well consult predicted or historical congestion, but they cannot observe other users' choices at the current instant, and so cannot anticipate the momentary fluctuations these choices produce. 
They are therefore unresponsive to short-term flow fluctuations while remaining systematically responsive to price signals—precisely the lever a pricing mechanism can exploit.

\paragraph{Decision-Dependency and Performativity}
Shifting to A1 changes the nature of the price–demand relationship: 
users no longer converge to an anticipatable equilibrium the platform can optimize over, and the environment becomes genuinely \textit{Decision-Dependent} (DD) \cite{perdomo2020performative}.

In classical stochastic optimization, uncertainty is exogenous, drawn from a fixed distribution that the optimizer cannot influence.
In our data-driven pricing context, however, the flow distribution is endogenous. When the platform sets a price $p$ (the decision variable), it does not merely select a point on a static demand curve but fundamentally alters the underlying distribution $\cD(p)$ from which network flows (the performative variable) are realized.
This phenomenon, Performativity \cite{perdomo2020performative}, means that the prediction itself shifts the distribution it aims to predict.
In our problem, the platform must learn to price against a shifting environment: each price adjustment intended to improve current flow patterns reshapes those patterns in response (see Fig.\,\ref{fig:feedback}, 
Fig.\,\ref{fig:DDDsterring}).
This feedback loop necessitates a re-evaluation of what ``optimal'' means.
For an AI system learning to price in such environments, the standard optimality concept, \textit{Performative Stability} fails catastrophically in our model.
We thus pursue \textit{Projected Performative Optimality}, a more rigorous standard that explicitly accounts for distributional shifts induced by the decision variables.

\paragraph{A New Problem Class}
We identify \textit{Projected Performative Optimization} ($\Pi$PO) as a generalization of standard Performative Optimization (PO) \cite{perdomo2020performative}, a paradigm for learning decisions in environments that react to those decisions: the distribution map itself contains a projection operator, $x \!\sim\!\Pi(\mathcal{D}(p))$.
If $\Pi\!\neq\!I$, new phenomena emerge (gating Jacobian, probabilistic smoothing, dead zone, and constraint qualification (CQ) failure,
see \cref{sec:algoframeworkSAA}-\cref{sec:algo_design}), which are structurally absent from PO and require new theory and algorithms\footnote{If $\Pi\!=\!I$, $\Pi$PO reduces to PO, and results in the standard literature apply.}.
Our Optimal Network Pricing (ONP) is the first concrete instance of $\Pi$PO; our analysis applies to any $\Pi$PO problem sharing this projected structure.

\paragraph{Contribution and Organization}
We formulate, analyze, and solve ONP under the Oblivious User Assumption as a learning problem within a DD Optimization framework.
Our contributions are as follows.
\begin{enumerate}[leftmargin=*]\setlength{\itemsep}{0pt}
\item \textbf{Theory.}
We introduce an ONP model (\cref{sec:ONP}) incorporating a DD Distribution (DDD) with a projection operator $\Pi$ that encodes capacity constraints.
The projection creates a nonsmooth, nonconvex optimization landscape with ``gating'' effects, where price signals are locally extinguished by capacity saturation; we handle this via subdifferentials and probabilistic smoothing.
In \cref{sec:POPS}, we define $\Pi$PO, derive its KKT conditions, and prove that PS collapses to a trivial solution in ONP.
We further prove the expected objective is strongly convex, guaranteeing a unique global $\Pi$PO solution despite the nonsmooth landscape.
Furthermore, we analyze the gap between $\Pi$PO and Approximate PO (APO) (\cref{sec:POPS:subsec:APO}), characterize regime-dependent sensitivity to the Dead Zone (DZ) and regularization parameter, and explain why the exact Jacobian is necessary.

\item \textbf{Algorithm.}
We propose a framework (\cref{sec:algoframeworkSAA}) that integrates SAA with TR-SQP to target $\Pi$PO rather than mere stability to learn $\Pi$PO prices from sampled demand responses.
\cref{sec:algoframeworkSAA} also establishes the gating Jacobian and provides sparse gradient/Hessian.
In \cref{sec:DZ_SampleComplexity}, we provide a sample complexity analysis, establishing theoretical bounds for approximation accuracy and ensuring  convergence of our framework.
The TR-SQP solver (\cref{sec:algo_design}) includes a DZ escape mechanism, necessary when CQs fail.

\item \textbf{Application.}
We show (\cref{sec:exp},\cref{sec:scalability}) that our $\Pi$PO-targeting algorithm, deployable as an AI decision-making system, significantly outperforms selected gradient-based methods, unlocking substantial social welfare gains by correctly accounting for gating effects (e.g., 
84\% on GÉANT). 

\end{enumerate}
We conclude in \cref{sec:conc} and summarize the paper in Fig.\,\ref{fig:organization}.
Our work draws on several disciplines\footnote{
Algorithmic mechanism design, 
concentration inequalities, 
convex analysis, 
decision-dependent modeling, graph theory, network pricing, stochastic optimization.
}, uses extensive notation and footnotes.
For readability, proofs and additional details are in the Appendix.

\paragraph{Basic notation}
``ONP'' and ``DD'' denote Optimal Network Pricing and Decision-Dependent, respectively; ``DDD'' means DD Distribution.
We work on a network (a simple strongly connected directed graph\footnote{No self-loops, no multi-edges, and a path exists between every pair of nodes.}) with vertex set $\cV$ (intersections, servers, hubs),
edge set $\cE$ (links connecting vertices: roads, cables, social interactions), and
route set $\cR \!\subset\! \textrm{PowerSet}(\cE)$ (paths through the network).
The assignment matrix $A \!\in\! \{0,1\}^{|\cR| \times |\cE|}$ encodes the route-edge relationship\footnote{$A_{re} \!=\! 1$ if edge $e \!\in\! \cE$ belongs to route $r \!\in\! \cR$, and $0$ otherwise.}.
This matrix is the fundamental topological operator of our framework, translating edge-wise congestion to route-wise costs\footnote{We express $\sum_{e \in r} (\cdot)$ as $\sum_{r \in \cR} A_{re} (\cdot)$ or $\sum_{e \in \cE} A_{re} (\cdot)$, depending on context.}.
Network flow demand is segmented into commodities $\cK$, each $k\!\in\!\cK$ representing an origin-destination pair.
The matrix $K\!\in\! \{0,1\}^{|\cK| \times |\cR|}$ maps routes to the commodities they serve\footnote{$K_{k,r}\!=\! 1$ if route $r$ connects the source and sink of commodity $k$.}.
The system has two variables: route flows $x \!\in\! \IR^{|\cR|}$ and platform-set route prices $p \!\in\! \IR^{|\cR|}$, bounded by $x \!\in\! [\zeros, x_{\max}]$, $p \!\in\! [p_{\min}, p_{\max}]$, and demand constraint $Kx \!\geq\! l$ ensuring total flow per commodity meets minimum demand $l$.
$\Pi_{[\zeros,x_{\max}]}(y)$ projects $y$ onto the box $[\zeros,x_{\max}]$.
Other notation is standard\footnote{$\langle a,b\rangle \!=\! a^\top b$ (inner product), $\|\cdot\|$ the $\ell_2$-norm. 
}.

\subsection{Literature Review}
Our work intersects behavior model, network pricing, and DD optimization.

\paragraph{Network User Behavior Modeling and Network Pricing}
Network user behavior modeling falls into three paradigms: equilibrium \cite{wardrop1952road}, utility \cite{hayrapetyan2005network}, and demand \cite{narang2023multiplayer}.
Classical approaches rely on A0 and formulate congestion games \cite{rosenthal1973class} where flows settle into equilibrium.\footnote{A state where no agent can unilaterally reduce their cost by switching routes.} While theoretically robust, these models focus on finding an equilibrium for fixed total demand.
Our work diverges by adopting the Oblivious User assumption (A1).
While \cite{karakostas2011degradation} uses A1 to study network performance degradation (Price of Anarchy \cite{roughgarden2002bad,roughgarden2005selfish}), they do not address the \textit{inverse problem}: optimizing pricing to shape flow distributions and steer oblivious users toward a system-optimal state.
We fill this gap by treating oblivious behavior as a predictable response function  for optimization within the DD framework.

Network pricing literature has mainly focused on Pigouvian taxes: tolls that internalize congestion externalities to steer systems toward a social optimum.\footnote{Algorithmic game theory often models these as Stackelberg games \cite{briest2012stackelberg, roch2005approximation}.}
Works such as \cite{hayrapetyan2005network} frame this as a mechanism design using disutility functions, but are typically limited to convex/linear formulations \cite{falkner2000overview, shakkottai2006economics} with simplified topologies and deterministic demand.\footnote{Works using parallel link topologies (\cite{hayrapetyan2005network}, \cite[Sec.3]{karakostas2011degradation}) or deterministic demand rarely integrate demand stochasticity with DDD constraints.}
We relax the parallel-link assumption to general graph topologies capturing complex route overlaps, and explicitly integrate high-dimensional stochastic demand.

\paragraph{DD Optimization}
Our work lies within DD optimization \cite{perdomo2020performative}, where flow distribution depends on the decision variable (price).
\cite{perdomo2020performative} formalized the distinction between \textit{Performative Stability} (PS) and \textit{Performative Optimality} (PO).
Most DD literature targets PS due to its tractability.\footnote{See \cite{mendler2020stochastic,  drusvyatskiy2023stochastic} for PS-targeting stochastic algorithms  under strong convexity, 
\cite{narang2023multiplayer} for multi-player extensions,
\cite{li2024stochastic} for relaxed convexity, and 
\cite{wang2025projected} for constraints on performative variables.}
However, we show PS is myopic in ONP, yielding a trivially suboptimal solution, and  PS-style updates fail structurally, motivating our $\Pi$PO-targeting algorithm.

We pursue $\Pi$PO\footnote{\cite{miller2021outside}, \cite{ray2022decision} also target PO, but on unconstrained smooth problems.
Our ONP is constrained, nonsmooth and nonconvex.} using machinery from \cref{sec:algoframeworkSAA}--\cref{sec:algo_design}, drawing on tools from several areas.
Our ONP presents fundamental difficulties unaddressed in the literature (Table~\ref{tab:DD_lit}): prior methods are inapplicable because most works do not consider hard constraints on the performative variable.
Moreover, projection inside the distribution map is unaddressed in the literature, including \cite{wang2025projected},\cite{ray2022decision}.
Our work thus contributes not only a new application but also a new \emph{problem class}, with direct implications for AI systems operating under performative, capacity-constrained feedback: $\Pi$PO subsumes all cited PO formulations as the special case $\Pi\!=\!I$, and probabilistic smoothing, the DZ geometry, and DZ escape mechanism are tools specific to $\Pi\! \neq\! I$ that have no counterpart in the existing PO literature.

\begin{table}[h!]
\centering
\caption{DD works.
C: constrained.
NcvxObj: nonconvex objective w.r.t. decision variable.
\cite{wood2021online},\cite{drusvyatskiy2023stochastic} handle constraints w.r.t. the decision variable but not the performative variable.
}
\label{tab:DD_lit}
\footnotesize
\begin{tabular}{l|l p{1.3cm} | p{1.85cm} p{4.2cm}}
\hline
\textbf{Work} & \textbf{C} & \textbf{NcvxObj} & \textbf{Optimality} & \textbf{Algorithm}
\\  \hline

\cite{perdomo2020performative},\cite{mendler2020stochastic},\cite{narang2023multiplayer},\cite{li2024stochastic} &  \no & \yes & PS & Retraining, SGD, GD \\ 

\cite{wood2021online},\cite{drusvyatskiy2023stochastic},\cite{wang2025projected} & 
\yes & \yes & 
PS &
Online, SGD, projected GD \\ 

\cite{miller2021outside} & 
\no & \no & 
approx.~PO  & 2-stage map estimation \\ 

\cite{ray2022decision} & 
\no & \yes & 
PO & epoch‑based 0th/1st‑order SGD \\ 

\hline
\textbf{Ours} & 
\yes & \yes & 
\textbf{Global $\bPi$PO} & SAA \& TR-SQP (2nd-order)
\end{tabular}
\end{table}
\section{Formulating Optimal Network Pricing (ONP)}\label{sec:ONP}
Our first contribution is the ONP problem, a principled formulation for an intelligent pricing system learning to allocate network resources under stochastic, capacity-constrained demand.
\begin{tcolorbox}
\vspace{-3mm}
\begin{equation}\label{prob:ONP}
\hspace{-4mm}
\begin{array}{cll}
\displaystyle \min_{p\in [p_{\min}, p_{\max}]} \hspace{-3mm} \IE_{x \sim \cD(p)}\!
\big[
\langle\tfrac{1}{2}Qx\!-\!s,x\rangle \!+\! \tfrac{\lambda}{2}\|p\|_2^2
\big]
~~\st~~
x\!\in\![\zeros, x_{\max}], ~ K x\!\geq\! l
\end{array}
\tag{ONP}
\end{equation}
\vspace{-6mm}
\begin{equation}\label{def:DDD}
~~~\cD(p) \!=\! \Pi_{[\zeros, x_{\max}]}( Bp \!+\! \zeta ), \qquad \zeta \sim \cN(\mu,\Sigma).\qquad 
\tag{DDD} 
\end{equation}
\end{tcolorbox}
\noindent
ONP and DDD form the foundation of this work.
Below we derive them from first principles, building from microscopic edge-wise congestion to macroscopic network flow distributions.
The derivation provides insights for later sections.

\subsection{Deriving the Quadratic Congestion Cost}\label{sec:ONP:subsec:derive}
Unlike abstract utility models with generic costs, we build the network-wide cost from linear edge-wise congestion{\footnote{This formulation is well-grounded in congestion games \cite{vickrey1969congestion, monderer1996potential, koutsoupias1999worst}.
Our methodology extends to nonlinear edge costs that preserve the convexity of the objective in \cref{thm:fE_strong_cvx}.
}.
Let $c^{\text{coe}}_e$ be the marginal cost coefficient\footnote{Marginal cost increase per unit additional flow.}
and $c^{\text{os}}_e$ the fixed offset\footnote{Cost at zero flow, e.g.\ free-flow travel time.}; 
the cost\footnote{E.g.\ latency, travel time, or fuel consumption.} experienced per unit flow on edge $e$ is
$\text{cost}_e(x)
\!=\! c^{\text{coe}}_e(\sum_{r \ni e}\!x_r)+\! c^{\text{os}}_e$.
Using $A$:
$\text{cost}_e(x)
\!=\!c^{\text{coe}}_e(\sum_{r\in\cR}\! A_{re} x_r)\!+\! c^{\text{os}}_e$.
User route cost is the sum of edge costs: $\text{cost}_r(x)\!=\! \sum_{e \in r}\!\text{cost}_e(x)$, so
\[
\text{cost}_r(x) 
\!=\! \sum_{e \in r}c^{\text{coe}}_e \hspace{-1mm}\left(\sum_{r' \in \cR} A_{r'e} x_{r'}\right) \!+\! c^{\text{os}}_e
\!=\! \sum_{e \in \cE} A_{re} \hspace{-1mm}\left( c^{\text{coe}}_e \sum_{r' \in \cR} A_{r'e} x_{r'} \!+\! c^{\text{os}}_e \right).
\]
From a mechanism design perspective, the platform minimizes the \textbf{total system cost} (social welfare) by summing costs across all flow units and routes.
\begin{proposition}\label{prop:CA}
The total system cost $\text{cost}_{\text{tot}}(x)\!=\!\sum_{r\in \cR}x_r \text{cost}_r (x)$ is the quadratic form
$\tfrac{1}{2} \langle Qx\!-\!s, ~x \rangle$
with $Q \!=\! 2A \Diag(c^{\text{coe}})A^\top \succeq \zeros$ and $s \!=\!-A c^{\text{os}}$.
\end{proposition}
\paragraph{Meaning of $Q$}
The symmetric PSD matrix $Q$ encodes \textit{Weighted Congestion Interaction}.
Since $(AA^\top)_{ij}\!=\!\sum_e A_{ie}A_{je}$, the off-diagonal entry $(AA^\top)_{ij}$ counts shared edges between routes $i$ and $j$, thus measuring their coupling.\footnote{$(AA^\top)_{ij}\!=\!0$ means routes $i,j$ are disjoint; $(AA^\top)_{ij}\!>\!0$ means flow on route $j$ contributes to congestion on route $i$.}
The diagonal entry $(AA^\top)_{rr}$ is the route length (number of edges).
Scaling by $\Diag(c^{\text{coe}})$ weights overlaps by cost sensitivity of shared edges.\footnote{A large $Q_{ij}$ means routes $i,j$ share edges where high flow on one significantly increases cost on the other.}
Our $Q$ generalizes parallel-link models \cite{hayrapetyan2005network,karakostas2011degradation}
\footnote{In these models, routes share no edges, simplifying $Q$ to a diagonal matrix.} to arbitrary graphs with complex path overlaps\footnote{
\begin{tikzpicture}[
 ->, thick, 
 node distance=0cm,
 every node/.style={circle, draw,inner sep=0.5pt},
 edge/.style={draw}
]
\node (s) {s};
\node (a) [right of=s, xshift=0.38cm] {1};
\node (c) [below right of=a, xshift=0.25cm, yshift=-0.32cm] {3};
\node (d) [right of=a, xshift=0.49cm] {2};
\node (t) [right of=d, xshift=0.39cm] {t};
\draw (s) -- (a);
\draw (a) -- (d);
\draw (a) -- (c);
\draw (c) -- (d);
\draw (d) -- (t);
\end{tikzpicture}
has paths $s$12$t$ and $s$132$t$ sharing edges $s$1 and 2$t$.
Parallel-edge models cannot represent such structures; ONP accommodates them naturally via $Q$.
}.

\subsection{The Decision-Dependent Distributional Demand}\label{sec:ONP:subsec:DDDD}
We now formalize the Oblivious User (A1): flow $x$ is not a strategic variable minimizing user cost but a stochastic variable drawn from a price-dependent distribution.
We propose \eqref{def:DDD} as our projected linear elasticity model, integrating price elasticity, demand uncertainty, and capacity constraints.

\paragraph{Price Elasticity}
The matrix $B \!\in\! \IR^{|\cR| \times |\cR|}$ encodes demand sensitivity to price, capturing oblivious demand under A1.
We assume the law of demand: increasing a route's price reduces its demand \cite{narang2023multiplayer}.
Mathematically, $B$ has a negative diagonal entries.
Off-diagonal entries $B_{ij}$ encode cross-route substitution.\footnote{If route $j$ substitutes for route $i$, a price increase on $j$ may raise demand for $i$ ($B_{ij}\!>\!0$).
For physical realism, we set $B_{ij} \!>\! 0$ only when routes $i$ and $j$ serve the same source-sink pair; otherwise $B_{ij} \!=\! 0$.}
$B$ maps price $p$ to \textit{deterministic latent demand} $Bp$, then corrupted by noise $\zeta$ to form \textit{latent demand} $y \!=\! Bp \!+\! \zeta$.

\paragraph{Demand Uncertainty}
$\zeta \!\sim\! \cN(\mu, \Sigma)$ models stochastic demand\footnote{$\zeta$ encodes unexplained factors (user habits, irrational decisions, real-time conditions).
We assume $\cN(\mu,\Sigma)$ for tractability; the framework extends to other continuous distributions.} and is independent of price.
The mean $\mu$ captures habitual baseline demand; $\Sigma$ models variance across users and scenarios.
We set $\Sigma \!\coloneqq\! \sigma^2 (I\!+\! \rho AA^\top)$ where $\sigma$ quantifies uncertainty magnitude, $AA^\top$ encodes topological route overlap, and $\rho$ ensures positive correlation across adjacent routes.

\paragraph{Projection's gating}
The projection $\Pi$ makes the mapping from $p$ to $x$ nonlinear, nonsmooth and  nonconvex:
this mapping is piecewise linear with kinks at the boundaries.
If $y \!\in\! (\zeros, x_{\max})$, then $\nabla_p x \!=\! B$, giving a linear price mechanism.
If $y \!\geq\! x_{\max}$, then $x \!=\! x_{\max}$ and $\nabla_p x \!=\! \zeros$: the price mechanism is gated off.\footnote{When latent demand exceeds $x_{\max}$, price changes do not affect flow.}
This poses a challenge (\cref{sec:POPS}), motivating our specialized algorithm (\cref{sec:algo_design}).

Both the projection and randomness are critical.
The projection encodes realistic flow bounds affecting pricing decisions:\footnote{
Without the upper bound, the platform may post artificially high prices to curb congestion, inflating the planner's objective.
Omitting the lower bound can force prices below the true optimum (due to regularization), increasing congestion cost. \label{fn:gating_effects}
}
it enforces operationally meaningful flow ranges and prevents extreme pricing that degrades the true objective.
Replacing $\zeta$ by its expectation is generally invalid, as $\IE[\Pi_{[\zeros,x_{\max}]}(y)]$ has no closed form:\footnote{$\Pi_{[\zeros,x_{\max}]}$ is nonlinear; computing its expectation requires moments of projected multivariate distributions involving multivariate CDFs.} projection does not commute with expectation\footnote{$\IE[\Pi(Bp\!+\!\zeta)]\neq \Pi(Bp\!+\!\IE[\zeta])$.
Projection truncates boundary mass and shifts the induced distribution; using $\IE[\zeta]$ ignores this truncation bias.} and a Jensen gap may appear.\footnote{For convex $\phi$, Jensen's inequality gives $\IE[\phi(X)] \!\geq\! \phi(\IE[X])$.}
We thus approximate $\IE[\Pi_{[\zeros,x_{\max}]}(y)]$ via Monte Carlo (\cref{sec:algoframeworkSAA:subsec:stearing}) to capture both projection nonlinearity and near-bound variability.

\paragraph{ONP problem}
Combining the cost and demand models gives \eqref{prob:ONP}: 
the platform minimizes expected total system cost (social welfare) plus price regularization (to prevent price gouging), subject to flow and demand constraints.
ONP is inherently a 2-agent problem: the \emph{platform} sets prices $p$, defines the objective, and enforces feasibility; the \emph{users} realize flows via the performative map $\cD(p)$.
As a \textbf{Constrained DD QP}, the expectation $\IE$ is over $\zeta$ indirectly via $x\!\sim\! \cD(p)$, and  $K\IE[x]\!\geq\! l$ is enforced in expectation.
We adopt assumption A2:
\begin{tcolorbox}
\textbf{(A2) Aligned bounds.}
The user-side and platform-side flow bounds coincide: $x \in [\zeros, x_{\max}]$ in both \eqref{prob:ONP} and \eqref{def:DDD}.
\end{tcolorbox}
\noindent
In general, the platform may enforce $x \!\in\! [\zeros, x_{\max}^{\text{pltf}}]$ while users operate under $\cD(p) \!=\! \Pi_{[\zeros, x_{\max}^{\text{usr}}]}(Bp\!+\!\zeta)$ with $x_{\max}^{\text{pltf}} \!\neq\! x_{\max}^{\text{usr}}$.
Relaxing A2 introduces a case split in the Dead Zone geometry (\cref{sec:POPS:subsec:Jac_DZ}) but does not alter the core analysis.\footnote{Under misaligned bounds, the DZ is governed by $x_{\max}^{\text{usr}}$ while feasibility is governed by $x_{\max}^{\text{pltf}}$. All results extend with notational overhead.}

Fig.\,\ref{fig:feedback} summarizes the ONP feedback process, a characteristic of AI systems acting in performative environments, where decisions reshape the very distribution being optimized.
\begin{figure}[h!]
\includegraphics[width=0.55\linewidth]{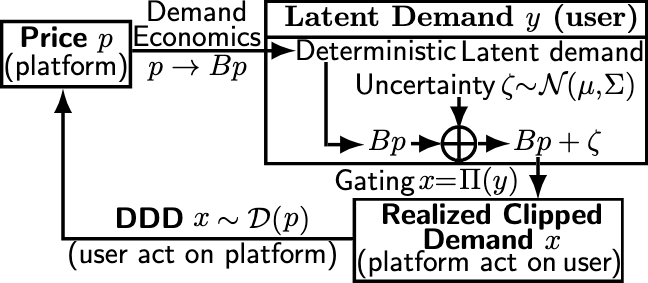}
\setlength{\tabcolsep}{-5pt}
\renewcommand{\arraystretch}{0.9}
\raisebox{0.85\height}{
$\hspace{-6pt}\begin{array}{l}
p : \text{price}
\\[-1pt]
B : \text{elasticity economics}
\\[-1pt]
Bp : \text{deterministic demand}
\\[-1pt]
y \!=\! Bp \!+\! \zeta : \text{latent flow}
\\[-1pt]
\zeta : \text{uncertainty}
\\[-1pt]
\mu, \Sigma : \text{mean and variance}
\\[-1pt]
x \!=\! \Pi_{[\zeros,x_{\max}]}(y) : \text{realized flow}
\\[-1pt]
\Pi : \text{projection / gating}
\end{array}$
}
\caption{The feedback mechanism relating the three variables: price $p$ (platform), latent demand flow $y$ (user), and realized flow $x$.
}
\label{fig:feedback}
\end{figure}
\section{SAA Formulation and Derivatives}\label{sec:algoframeworkSAA}
We establish the computational foundations for solving ONP as a data-driven learning problem.
The objective is an expectation over a projected distribution without a closed-form solution, and the embedded projection induces a nonsmooth, nonconvex landscape.
We address these challenges in three steps:
\cref{sec:algoframeworkSAA:subsec:stearing} interprets pricing as ``flow steering'' and applies Sample Average Approximation (SAA) to discretize the stochastic objective;
\cref{sec:algoframeworkSAA:subsec:JacobianSmoothing} introduces the Gating Jacobian and proves almost sure smoothing;
\cref{sec:DZ_SampleComplexity:subsec:sparse} derives gradients and Hessians, exploiting sparsity for scalability.
The resulting algorithm is presented in \cref{sec:algo_design}.

\subsection{Pricing as Flow Steering and SAA}\label{sec:algoframeworkSAA:subsec:stearing}
ONP treats price $p$ as both a variable to optimize and a steering mechanism: 
by changing $p$, the platform shifts latent demand $y$ so that realized flows $x$ minimize congestion while respecting constraints (see Fig.\,\ref{fig:feedback}, 
Fig.\,\ref{fig:DDDsterring}).
Mathematically, we substitute $x$ by the performative map $x_{p,\zeta} \!=\! \Pi_{[\zeros, x_{\max}]}(Bp\!+\!\zeta)$,
transforming \eqref{prob:ONP} from optimization over $(x,p)$ into a problem over $p$ alone\footnote{The substitution removes the constraint $x \in [\zeros, x_{\max}]$ from ONP due to (A2), shifting the task from resource allocation to distributional control.}:
\begin{equation}\label{prob:ONP-nox}
\underset{p\in [p_{\min}, p_{\max}]}{\argmin}
 \ftrue \coloneqq 
\IE [\langle \tfrac{1}{2}Qx_{p,\zeta} \!-\! s, ~ x_{p,\zeta} \rangle]
\!+\! \tfrac{\lambda}{2}\|p\|_2^2
~~ \st\ K \IE[x_{p,\zeta}] \geq l.
\tag{ONP$^{-x}_{\IE}$}
\end{equation}
\noindent
This ``steering'' capability distinguishes $\Pi$PO (\cref{sec:POPS}) from PS-seeking methods, as our algorithm explicitly accounts for how demand responds to price changes.

The expectation $\IE$ in \eqref{prob:ONP-nox} has no closed form.
We approximate it via SAA \cite{shapiro2021lectures}: sample $N$ i.i.d.\ scenarios $\zeta^1,\ldots,\zeta^N$.
Denoting $\Avg[y^i] \!\coloneqq\! \tfrac{1}{N}\sum_i y^i$ as the empirical average,
we have
\begin{equation}\label{prob:ONP-nox-SAA}
\hspace{-1mm}
\textstyle
\underset{p\in [p_{\min}, p_{\max}]}{\argmin}
f_N(p) 
\hspace{-0.5mm} = \hspace{-0.5mm} \Avg\!
[\langle \tfrac{1}{2}Q x^i_p \!-\! s, ~ x^i_p \rangle]
\!+\! \tfrac{\lambda}{2} \|p\|_2^2
~~\st~~ c_N(p) \leq \zeros,
\tag{ONP$_N$}
\end{equation}
where $y^i\!=\!Bp\!+\!\zeta^i$, $x^i_p\!=\! \Pi_{[\zeros, x_{\max}]}(y^i)$ is the realized flow in scenario $i$, $f_N$ is the SAA objective, and $c_N(p) \!\coloneqq\! l \!-\! K \Avg[x^i_p]$ is the SAA constraint.
SAA is grounded in its asymptotic properties (\cref{sec:DZ_SampleComplexity}).\footnote{
We do not resample $\zeta^i$ at each iteration:
\textbf{1.} TR-SQP (\cref{sec:algo_design}) require the \emph{same} $f_N$ at the outset.
\textbf{2.} Per-step resampling reduces to RRM \cite{perdomo2020performative},
which fails on ONP (Theorem~\ref{thm:PS_failure}).
}
\eqref{prob:ONP-nox-SAA} is NP-hard to solve exactly\footnote{Optimization of constrained, nonconvex piecewise-quadratic programs is NP-hard \cite{sahni1974computationally,pardalos1991quadratic}.},
but has polynomial memory and sample complexity (\cref{sec:DZ_SampleComplexity}, \cref{sec:app:subsec:sparse}).

\subsection{The ``Gating'' Jacobian and Almost Sure Smoothing}\label{sec:algoframeworkSAA:subsec:JacobianSmoothing}
To minimize $f_N(p)$, we differentiate through $\Pi_{[\zeros,x_{\max}]}$ via subdifferential \cite{clarke1975generalized}.
\begin{definition}
The diagonal generalized Jacobian $J \!=\! \nabla_p \Pi_{[\zeros, x_{\max}]}(y)$ is 
\begin{equation}\label{eq:Jac_subdiff}
J_{rr}
=
\begin{cases}
\text{singleton } \{0\} & y_r \notin [0,(x_{\max})_r] \quad (\text{Inactive/Saturated}), \\[-1.5pt]
\text{singleton }  \{1\} & y_r \in (0,(x_{\max})_r) \quad (\text{Active/Linear}), \\[-1.5pt]
\text{interval } [0,1] & y_r \in \{0, (x_{\max})_r\} \quad (\text{Boundary/Kink}).
\end{cases}
\end{equation}
\end{definition}
\noindent
Nondifferentiability at boundary points $\{\zeros, x_{\max}\}$ introduces kinks that complicate gradient analysis \cite{rockafellar1970convex,rockafellar1998variational}.
We resolve this by \textit{smoothing}:\footnote{The Gaussian density of $\zeta$ makes the probability of hitting a kink zero.
In finite-precision arithmetic this probability is ${\sim}10^{-16}$, negligible.}
\begin{lemma}[Probabilistic Smoothing]\label{lem:as_smoothing}
For $\zeta \!\sim\! \cN(\mu, \Sigma)$ with $\Sigma \!\succ\! \zeros$, the set where the derivative is undefined has measure zero: $\IP( y_r \in \{0, (x_{\max})_r\}) \!=\! 0$.
\end{lemma}
\noindent
This converts a combinatorially nonsmooth problem into smooth stochastic optimization: both $\ftrue$ and $f_N$ are differentiable \textit{almost everywhere}.
Since $J$ is diagonal with binary entries almost surely, we use Random Matrix Theory to analyze sample complexity (\cref{sec:DZ_SampleComplexity}).

\subsection{Gradient, Hessian and Efficient Computation}\label{sec:DZ_SampleComplexity:subsec:sparse}
With probabilistic smoothing established, we compute the analytical gradient and Hessian of $f_N$ and $c_N$ almost surely.
\begin{lemma}\label{gradientHessian_derivation}
$\nabla_p f_N(p) \!=\! B^\top \Avg[J^i (Qx^i_p \!-\! s)] \!+\! \lambda p $ where $J^i$ is the Jacobian for sample $i$, and $\nabla_p^2 f_N(p) \!\overset{a.s.}{=}\! B^\top \Avg[J^i Q J^i] B \!+\! \lambda I$.
Also 
$\nabla_p c_N(p) \!=\! \!-\! K \Avg[J^i] B$ and $\nabla_p^2 c_N(p) \!\overset{\text{a.s.}}{=}\! \zeros$.
\end{lemma}
\noindent
The marginal congestion cost $Qx^i_p\!-\!s$ for scenario $i$ is propagated back to prices by $B$.
Crucially, $J^i$ acts as a \textbf{gate}: if route $r$ is saturated in scenario $i$, then $J^i_{rr} \!=\! 0$ and the congestion signal from $r$ is lost.
The gating $J$ characterizes the performative effect (the endogenous feedback induced by the current decision $p$) and is the key to Performative Optimality (\cref{sec:POPS}).
This environmental shift recursively influences the optimization of $p$ in subsequent iterations (see Fig.\,\ref{fig:DDDsterring}).

\paragraph{Sparsity-Aware Computation}
When computing $Q \!=\! 2 A \text{Diag}(c^{\text{coe}}) A^\top$, we use a matrix-free product exploiting the sparsity of $A$\footnote{Each row of $A$ has nonzero entries only for edges it traverses.} (\cref{sec:app:subsec:sparse}).

\section{Sample Complexity and Dead Zone}\label{sec:DZ_SampleComplexity}
The projection $\Pi_{[\zeros,x_{\max}]}$ in the DDD map introduces \emph{Dead Zones} (where the SAA gradient loses all congestion information) in the price space.
We now quantify sample complexity, define the Dead Zone, prove that no minimizer lies inside it, and derive the sample thresholds that govern algorithm design in \cref{sec:algo_design}.
These results inform how much data an intelligent pricing system needs to learn reliable prices.

\paragraph{SAA Sample Complexity}\label{sec:DZ_SampleComplexity:subsec:SAA}
Moving from \eqref{prob:ONP-nox} to \eqref{prob:ONP-nox-SAA} introduces approximation error.
We quantify how $N$ must scale through three levels:
\textbf{Consistency} ($N_{SAA}$): ensure $f_N$ approximates $\ftrue$ and their differentials correctly.
\textbf{Existence} ($N_{\DZ}$): prevent gradient vanishing due to \DZ.
\textbf{Stability} ($N_{\text{op}}$): guarantee TR-SQP convergence.
We start with the accuracy of the SAA. 

\begin{proposition}[$N_{SAA}$ by Uniform Convergence and Error Gap]\label{prop:concentration_combined}
Let $\cP$ be the box constraint of $p$.
For ONP, w.p.  at least $1\!-\!\delta$, we have that
\[
\sup_{p \in \cP} |f_N - \ftrue| \leq \epsilon,
~~~
\sup_{p \in \cP} \|\nabla f_N - \nabla \ftrue\|_2 \leq \epsilon,
~~~
\sup_{p \in \cP} \|\nabla^2 f_N - \nabla^2 \ftrue\|_2 \leq \epsilon,
\] 
with a constant $L$ (see Appendix), it suffices to take
\begin{tcolorbox}\vspace{0mm}
\[
N_{SAA}
\!=\! \cO\Big(\epsilon^{-2}\big[
 \big(\ln\tfrac{\|p_{\max}-p_{\min}\|_2 L}{\epsilon}\big) |\cR| + \ln\tfrac{|\cR|}{\delta}
\big]
\Big)
\tag{$N_{SAA}$}
\]
\end{tcolorbox}
\noindent
\textbf{Error Gap}: $\displaystyle |\min_{p \in \cP} f_N(p) \!-\! \min_{p \in \cP} \ftrue(p)| \!\leq\! \epsilon$
w.p. $1\!-\!\delta$ if $N \!\geq\! N_{SAA}$.
\end{proposition}
\noindent
Prop.~\ref{prop:concentration_combined} \textit{transfers} property between $\ftrue$ and $f_N$.
Function bound gives the gap $|\min\! f_N \!-\!\min\!\ftrue|\!\le\!\epsilon$.
The gradient bound enables stationarity transfer to prove Corollary~\ref{cor:end_to_end}: if $\|\nabla f_N(p_N^*)\|\!\le\!\epsilon_{\mathrm{opt}}$, then $\|\nabla\ftrue(p_N^*)\|\!\le\!\epsilon_{\mathrm{opt}}+\epsilon$.
The Hessian bound enables curvature transfer for efficient algorithm execution (\cref{sec:DZ_SampleComplexity:subsec:Nop}).

\subsection{Route Activation and Dead Zone}\label{sec:POPS:subsec:Jac_DZ}
We now move to $N_{\DZ}$, where we define Route Activation Probability.
\begin{definition}[Route Activation Probability]\label{def:route_activation}
$\pi_r\coloneqq\IE[J_{rr}]$.    
\end{definition}
\noindent 
In deep congestion regimes, $\pi_r$ is small. 
For small $N$, the entries $J_{rr}^i$ may vanish across all samples; congestion information is lost and the solver stagnates.
Treating route activation as a Bernoulli event, the probability that route $r$ is dead (inactive) across all $N$ samples is $(1\!-\!\pi_r)^N$.
To observe at least one active scenario with confidence $1\!-\!\delta$, it suffices to take\footnote{From $1-(1\!-\!\pi_r)^N\ge 1\!-\!\delta$.}
\begin{tcolorbox}\vspace{-1mm}
\[
N_{\DZ} \geq |\ln(1-\pi_{\min})|^{-1}\ln\tfrac{1}{\delta}.
\tag{$N_{\DZ}$}
\]
\end{tcolorbox}
This is formalized by the \emph{Dead Zone}. 

\paragraph{Dead Zone}
The gradient $\nabla_p f_N$ (Lemma~\ref{gradientHessian_derivation}) decomposes into a regularization term $\lambda p$ and a stochastic congestion term $B^\top \Avg[J^i(Qx^i\!-\!s)]$.
The Jacobian \eqref{eq:Jac_subdiff} reveals the \textbf{Dead Zone} (\DZ): a \textit{random} region where latent demand is so extreme that $J^i\!=\!\zeros$.\footnote{$J_{rr}^i=0 \iff y_r^i\ge(x_{\max})_r$ or $y_r^i\le 0$.}
In this region, flow is completely insensitive to price; the congestion term vanishes, and $\nabla_p f_N$ collapses to $\lambda p$.
\begin{tcolorbox}\vspace{0mm}
\begin{equation}\label{def:DZ}
\begin{array}{rll}
\text{Partial} & \DZ_r(f_N) & \coloneqq
\{ p \in [p_{\min}, p_{\max}]  :  \Avg[J^i_{rr}(p, \zeta^i)] = 0 
\},
\\[2pt]
\text{Full} & \DZ(f_N) & \coloneqq
\{ p \in [p_{\min}, p_{\max}]  :  \Avg[J^i(p, \zeta^i)] = \zeros 
\}.
\end{array}
\tag{\DZ}
\end{equation}
\end{tcolorbox}

\paragraph{Interior and kink}
We partition \DZ\ into its open interior ($\mathrm{int}\DZ$) and its measure-zero boundary, the \emph{kink} ($\partial\DZ$).
The \emph{DZ-residual} $\mathrm{res}\coloneqq y-\Pi_{[\zeros,x_{\max}]}(y)$ quantifies displacement of latent demand from the feasible flow box.
\begin{itemize}[leftmargin=*]\setlength{\itemsep}{0pt}
\item \textbf{Interior ($p_r\in\mathrm{int}\DZ_r$):}  
$y_r \notin [0,(x_{\max})_r]$, so $J_{rr}\!=\!0$, $\mathrm{res}_r\!\neq\! 0$.
The gradient is blind to congestion, but $\mathrm{res}_r\!\neq\! 0$ provides a geometric restoring force.

\item \textbf{Kink ($p_r\in\partial\DZ_r$):}
$y_r$ sits exactly on the bound, so $\mathrm{res}_r\!=\!0$ and $J_{rr}\!\in\![0,1]$.
No restoring force exists; trust-region shrinkage handles this case (\cref{sec:algo_design}).
\end{itemize}
\noindent
A consequence of \eqref{def:DZ} and Lemma~\ref{gradientHessian_derivation}:
inside $\mathrm{int}\DZ(f_N)$, $J^i \!=\! \zeros$ for all $i$, so
\begin{equation}\label{eq:DZ_grad_collapse}
\nabla_p f_N(p) \!=\! \lambda p
\qquad \text{and} \qquad
\nabla_p c_N(p) \!=\! \zeros
\qquad\text{(in } \mathrm{int}\DZ\text{)}.
\end{equation}
The congestion signal vanishes entirely.

\paragraph{Probability of touching DZ}
As \DZ\ is a random set, we can talk about the chance of a point touching \DZ.
For small $\pi_{\min}$, a prohibitively huge number of samples may be needed to ensure that $p_N$ lies outside \DZ.
In particular, if $N\!<\!N_{\DZ}$, the probability that $p_N$ touches \DZ\ is substantial.
However, as we show below, the minimizer cannot lie in \DZ.

\begin{proposition}\label{prop:ftrue_no_DZ_min}
The minimizer of $f_N$, denoted as $p^*_N$, cannot lie in $\mathrm{int}\DZ(f_N)$. 
\end{proposition}
\noindent 
While $f_N$ has \DZ, $\ftrue$ has empty \DZ\ (\cref{sec:app:pf:prop:DZfE_no_PO_pt_pi_r_positive}). But note that for the global minimizer $p^*_{\IE}$ of $\ftrue$, Prop.~\ref{prop:ftrue_no_DZ_min} alone does not exclude $\IP\big( p^*_{\IE} \!\in\! \DZ(f_N) \big)\!=\!0$.
However, based on Prop.~\ref{prop:concentration_combined}, we have the following corollary:
\begin{corollary}\label{cor:pIE_not_in_DZ}
For $N \geq N_{SAA}(\epsilon, \delta)$, we have
$\IP\left(p^*_{\IE} \in \DZ(f_N)\right) \leq \delta.$
\end{corollary}
\begin{remark}[Trading Samples for Iterations]\label{rem:N1_algo}
Note that 
\begin{enumerate}[leftmargin=*]\setlength{\itemsep}{0pt}
\item Corollary~\ref{cor:pIE_not_in_DZ} tells $p^*_{\IE} \notin \mathrm{int}\DZ(f_N)$ w.h.p., but during optimization $p$ can get into $\mathrm{int}\DZ(f_N)$. 

\item $N_{\DZ}$ can be huge if $\pi_{\min}\approx 0$.
Fixing it requires exhaustive resampling.
\end{enumerate}
They motivate \textbf{DZ escape} (\textbf{reflection}) in Algorithm~\ref{algo:TRSQP} that \emph{trades sample complexity for
iteration cost}; see Lemma~\ref{lem:finite_dz_escape}.
\end{remark}

\subsection{Operational Sample Complexity}\label{sec:DZ_SampleComplexity:subsec:Nop}
TR-SQP operates on $f_N$, and requires positive curvature to make stable progress.
Since $\nabla^2 f_N(p) \overset{\text{a.s.}}{\succeq} \lambda I$ (Lemma~\ref{gradientHessian_derivation}), the curvature of $f_N$ is at least $\lambda$. Because $J^i$ is an unknown random matrix, $\lambda$ is the best deterministic lower bound on $\nabla_p^2 f_N(p)$.
Relying on $\lambda$ for curvature makes TR‑SQP arbitrarily slow for small $\lambda$ and disconnects convergence speed from the congestion structure.

The solution is to transfer curvature from $\ftrue$: we have that 
$\ftrue$ is $\alpha$-strongly convex
with $\alpha \!\propto\! \pi_{\min}^2$ (Theorem~\ref{thm:fE_strong_cvx}).
The Hessian bound in Prop.~\ref{prop:concentration_combined}  transfers this curvature to $f_N$: for $N \!\geq\! N_{SAA}(\epsilon,\delta)$,
we have $ \nabla^2 f_N(p) \!\succeq\! (\alpha \!-\! \epsilon)I$ almost surely.
If the approximation error smaller than the curvature ($\alpha \!\gg\! \epsilon$), then $f_N$ inherits the curvature of $\ftrue$.
Setting $\epsilon \!=\! \cO(\pi_{\min}^2) \!=\! \cO(\alpha)$ as the worst-case
binding constraint and substituting into $N_{SAA}$ gives a \textit{quartic} complexity as
\begin{tcolorbox}\vspace{0mm}
\[
N_{\text{op}}
~=~
\cO\left(\frac{1}{\pi_{\min}^4}
\left[|\mathcal{R}|\ln\frac{1}{\pi_{\min}^2} 
+ \ln\frac{|\mathcal{R}|}{\delta}\right]\right)
~\overset{\pi_{\min}\!\to 0}{=}~
\tilde{\cO}\left(\frac{1}{\pi_{\min}^4}\right).
\tag{$N_{\text{op}}$}
\]
\end{tcolorbox}
\noindent
$N_{\text{op}}$ tells how many samples are needed so that $f_N$ converges as fast as minimizing $\ftrue$ by  \textbf{recovering $\ftrue$'s
landscape} inside the SAA objective. 

\begin{remark}[$N_{\DZ}$ vs $N_{\text{op}}$]
In the regime $\pi_{\min}\!\to\! 0$ we have $|\ln(1-\pi_{\min})| \approx \pi_{\min}$, so $N_{\DZ}$ scales as $\cO(\pi_{\min}^{-1})$ and is dominated by $N_{\text{op}}$.
Note that $N_{\text{op}}$ guarantees stable progress toward the global optimum $p_{\IE}^*$, with convergence governed by $\pi_{\min}$ and $\alpha \propto \pi_{\min}^2$.
In ONP,  $\pi_{\min}$ is a key parameter: 
a large $\pi_{\min}$ reduces $N_{\text{op}}$ quartically while increasing curvature $\alpha$, yielding fast convergence; 
conversely, a small $\pi_{\min}$ makes curvature vanish, making convergence to $p_{\IE}^*$ difficult.
\end{remark}

\section{Optimality and the Failure of Stability}\label{sec:POPS}
In \eqref{prob:ONP} both $x$, $p$, and the distribution are endogenous, making ONP a \emph{decision-dependent} (DD) optimization problem, a core setting for AI systems that learn to act in performative environments.
In \cref{sec:POPS:subsec:KKT}, we first review the standard PS and PO notions from the DD optimization literature; then, using the dead-zone geometry in \cref{sec:DZ_SampleComplexity}, we introduce a new optimality concept called \ref{prob:general_proj_form} and derive its optimality conditions. 
We then introduce a baseline called Approximate PO (APO) in \cref{sec:POPS:subsec:APO},
 and analyze its gap to \ref{prob:general_proj_form}. 
Lastly, in \cref{sec:POPS:subsec:PS_collapse}, we show that PS degenerates in ONP,
and prove that PS and \ref{prob:general_proj_form} coincide under narrow conditions.

\subsection{Optimality Notions and KKT Conditions}\label{sec:POPS:subsec:KKT}
Two ways to define ``the solution'' in DD optimization arise from how the feedback loop between $p$ and $\cD(p)$ is treated \cite{perdomo2020performative}.
\begin{itemize}[leftmargin=*]\setlength{\itemsep}{0pt}
\item \textbf{Performative Stability (PS).}
  A PS price $p_{PS}$ is optimal for the fixed distribution itself induces:
  $p_{PS}=\mathrm{argmin}_p\IE_{x\sim\cD(p_{PS})}[\text{cost}_{\text{tot}}(x,p)]$.
  PS corresponds to iterative retraining strategies%
  \footnote{E.g.\ Repeated Risk Minimization \cite{perdomo2020performative}, Repeated Minimization \cite{drusvyatskiy2023stochastic}.}
  that treat the current distribution as static at each step, ignoring the full DD map.

\item \textbf{Performative Optimality (PO).}
  PO minimizes cost by explicitly accounting for how $\cD(p)$ shifts with $p$:
  $p_{PO}\!=\!\mathrm{argmin}_p \IE_{x\sim\cD(p)}[\text{cost}_{\text{tot}}(x,p)]$.
  Achieving PO requires differentiating through $\cD(p)$, anticipating the system's reaction to any price change.
\end{itemize}
The ONP objective is an instance of the projected performative problem
\begin{tcolorbox}\vspace{-1mm}
\begin{equation}\label{prob:general_proj_form}
\textstyle
p_{\Pi\text{PO}}
~=~\argmin_p\;\IE_{x \sim \Pi(\cD(p))} \left[\ell_1(x)\right] + \lambda \ell_2(p),
\tag{$\Pi$PO}
\end{equation}
\end{tcolorbox}
\noindent
where the performative map itself contains a projection operator.
Standard PO \cite{perdomo2020performative} is the special case $\Pi\!=\!I$, in which the distribution map is smooth, the Jacobian is never gated, no Dead Zone exists, and CQ holds everywhere.
The projection $\Pi \!\neq\! I$ breaks all four properties simultaneously, placing $\Pi$PO strictly outside the scope of existing PO theory
\cite{perdomo2020performative,mendler2020stochastic,miller2021outside,
wood2021online,drusvyatskiy2023stochastic,narang2023multiplayer,
li2024stochastic,wang2025projected}.

Below we move to the KKT conditions.
We focus on the interior case
($p_{\min} \!<\! p^* \!<\! p_{\max}$), omitting the box-constraint in $p$  for clarity.\footnote{Incorporating them adds complementarity
conditions $\mu_{\min} \!\geq\! \zeros, \mu_{\max} \!\geq\! \zeros$ with
$\mu_{\min} \odot p \!=\! \zeros, \mu_{\max} \odot (p_{\max} - p) \!=\! \zeros$
to the KKT system.
This lengthens expressions but does not alter argument, as the core analysis concerns the congestion term
$B^\top\IE[J(Qx-s)]$ and the DZ geometry, both independent of boundary activation.}

\paragraph{KKT conditions for $\mathit{\Pi}$PO (on $\ftrue$)}
As discussed after Prop.~\ref{prop:ftrue_no_DZ_min}, $\ftrue$ has no \DZ\ and CQs hold everywhere, so KKT conditions are necessary and sufficient for the unique global $\Pi$PO (Theorem~\ref{thm:fE_strong_cvx}).
Let $\gamma_{\Pi\text{PO}}$ be the Lagrange multipliers for demand; $p^*_{\Pi\text{PO}}$ satisfies
\begin{tcolorbox}\vspace{0mm}
\begin{equation}\label{eq:KKT_true}
B^\top \IE \bigl[J (Q x - s)\bigr]
- B^\top \IE[J] K^\top \gamma_{\Pi\text{PO}}
+ \lambda p^*_{\Pi\text{PO}} = \zeros.
\tag{$\Pi$PO-KKT$_{\IE}$}
\end{equation}
\end{tcolorbox}
\noindent
$\IE[J]\!=\!\Diag(\pi)$ gates congestion gradients on average:
pricing signals transmit only through routes with $\pi_r\!>\!0$, never through permanently saturated ones.
\paragraph{Clarke-stationarity on $f_N$}
For finite $N$, \DZ\ exist where sample Jacobians vanish entirely.
We analyze its effect on constraint violation.

\begin{proposition}[CQ Failure in partial DZ]\label{prop:fail_CQ}
If $p_r\in\DZ_r(f_N)$, both LICQ and MFCQ fail: $J^i_{rr}\!=\!0$ for route $r$ and all samples gives $[\nabla_p c_N]_r\!=\!0$ for route $r$.
\end{proposition}
\noindent 
For iterates in partial \DZ\ (some routes gated),
the affected constraint rows are handled by the
restoration step of Algorithm~\ref{algo:TRSQP}
(see \cref{sec:algo_design}).

In regions where $p_r \notin \DZ_r(f_N)$ for all $r$
(i.e.\ outside every partial \DZ),
CQ holds and Clarke-stationarity reduces to the SAA KKT:
\begin{tcolorbox}\vspace{-4mm}
\[
B^\top\Avg\bigl[J^i(Qx^i-s)\bigr]
\;-\; B^\top\Avg[J^i]\,K^\top\gamma_N
\;+\; \lambda p_N^*
\;\overset{\text{a.s.}}{=}\; \zeros.
\tag{$\Pi$PO-KKT$_N$}
\]
\end{tcolorbox}

\paragraph{Global $\mathit{\Pi}$PO of $\ftrue$}
While $f_N$ is nonconvex for finite $N$, $\ftrue$ is strongly convex, elevating local KKT stationarity to global optimality, based on an economically and mathematically natural assumption A3\footnote{
If $B$ is invertible, $\rk(KB) = \rk(K)$, so $K$ need only have full row rank. 
This standard condition guarantees no redundant/contradictory demand segments in $\cK$, and that $B$ preserves the linear independence. 
So pricing can influence every commodity flow, which is verifiable.
} :
\begin{tcolorbox}
\textbf{(A3) Demand responsiveness.}
The matrix $KB$ has full row rank:
each demand constraint responds to at least one price via the elasticity map $B$.
\end{tcolorbox}

\begin{theorem}[Probabilistic Strong Convexity]\label{thm:fE_strong_cvx}
Let $\pi_{\min}\!\coloneqq\!\min_{r\in\cR}\pi_r$ with $\pi_r\!=\!\IE[J_{rr}]$.
Assume $\pi_{\min}\!>\!0$ and $\lambda\!>\!0$.
Then $\nabla_p^2\ftrue(p)\!\succeq\!\alpha I$ with $\alpha\!=\!\lambda\!+\!\pi_{\min}^2\lambda_{\min}(Q)\sigma_{\min}^2(B)\!>\!0$.
As a result, $\ftrue$ is $\alpha$-strongly convex and admits a \textbf{unique global minimizer} $p^*_{\Pi\text{PO}}$, the unique solution to \eqref{eq:KKT_true}; there are no spurious local minima or saddle points.
\end{theorem}
\noindent
As $f_N\!\to\!\ftrue$ uniformly (Prop.~\ref{prop:concentration_combined}), solving \eqref{prob:ONP-nox-SAA} with sufficiently large $N$ yields a solution converging to $p^*_{\Pi\text{PO}}$ with high probability (Corollary~\ref{cor:end_to_end}).

\subsection{Approximate PO: applying PO methodology on $\mathit{\Pi}$PO}
\label{sec:POPS:subsec:APO}
We define a baseline that accounts for performative shifts but ignores capacity bounds by setting $J\!\approx\! I$.

\begin{definition}[Approximate PO (APO)]\label{def:approxPO}
The APO point solves the $\mathit{\Pi}$PO objective but computes gradients as if $\Pi \!=\! I$ (i.e., $\nabla_p x \!=\! B$, ignoring projection).
Its KKT condition is $B^\top\IE[Qx-s]-B^\top K^\top\gamma_{\text{APO}}\!+\!\lambda p_{\text{APO}}\!=\!\zeros$.
\end{definition}
\noindent
APO is equivalent to applying standard PO methodology, which assumes a smooth distribution map, to a $\Pi$PO problem.
It cannot distinguish price-responsive routes from gated ones, leading to systematic mispricing quantified by the gap analysis below.
\paragraph{The APO-$\mathit{\Pi}$PO Gap}
We quantify the inefficiency of APO via $\texttt{Gap} \!=\! \ftrue(p_{APO}) - \ftrue(p_{\Pi\text{PO}})$.
The gap arises because APO ignores the gating Jacobian: the flow $\hat{x}_{APO}$ \textit{perceived} by APO differs from the realized flow $x \!=\! \Pi_{[\zeros,x_{\max}]}(y)$ \textit{experienced} by users.
We assume (A4) for this analysis:
\begin{tcolorbox}
\textbf{(A4) Diagonal dominance.}
Own-route congestion-pricing feedback dominates cross-route substitution on shared edges: $-B^\top Q \!\geq\! \zeros$ element-wise.\footnote{A joint condition on $B$ and $Q$, distinct from $B\!\prec\! \zeros$, verifiable for any given network.
It holds when $|B_{ii}|$ is sufficiently large relative to $\sum_{j\neq i}|B_{ij}|$ and the structure in $Q$.}
\end{tcolorbox}
\noindent The gap exhibits regime-dependent sensitivity to $\lambda$ (Prop.~\ref{prop:reg_amp_gap}). 
In the \textbf{upper saturation regime} (small $x_{\max}$), APO overestimates demand and overprices; the quadratic price penalty amplifies \texttt{Gap} as $\lambda$ increases.
In the \textbf{lower saturation regime} (large $x_{\min}$), APO underestimates demand and underprices; congestion cost amplifies \texttt{Gap} as $\lambda$ decreases.

\begin{proposition}[Regime-dependent Gap]\label{prop:reg_amp_gap}
Under A4: (i) Upper regime: APO overflows and overprices ($\hat{x}_{APO}\!\geq\!x_{\mathit{\Pi}PO}$, $p^*_{APO}\!\geq\!p^*_{\mathit{\Pi}PO}$); larger $\lambda$ widens \texttt{Gap}.
(ii) Lower regime: APO underflows ($\hat{x}_{APO}\!\leq\!x_{\mathit{\Pi}PO}$); smaller $\lambda$ widens \texttt{Gap}.
\end{proposition}
\noindent
At the critical $\lambda^*$ balancing both regimes, \texttt{Gap} purely reflects the performative feedback, and we have the following lemma:
\begin{lemma}[Gradient mismatch]\label{lem:grad_mismatch}
$\|\nabla_p \ftrue(p_{APO})\| \!=\! \|B^\top \IE[(J - I)(Qx- s)]\|$, and
$\texttt{Gap} \!=\! \tfrac{1}{2\alpha} \|\nabla \ftrue(p_{APO})\|^2$ with
$\alpha \in [\lambda_{\min}(\nabla_p^2 \ftrue(p_{\Pi\text{PO}})),\, \lambda_{\max}(\nabla_p^2 \ftrue(p_{\Pi\text{PO}}))]$.
\end{lemma}

\begin{proposition}[Gap Amplification]\label{thm:global_gap_amp}
At $\lambda^*$, decreasing $\pi_r$ widens \texttt{Gap}: as $\pi_r$ drops, $(I\!-\!J)\to I$ amplifies the numerator $\|\nabla\ftrue(p_{APO})\|^2$ while the curvature $\alpha$ in the denominator decreases (Theorem~\ref{thm:fE_strong_cvx}), compounding the gap.
\end{proposition}
\noindent
These results justify the exact Jacobian treatment in our TR-SQP solver: ignoring gating leads to systematic mispricing across all $\lambda$ regimes, with the distortion worsening precisely when flow is most gated ($\pi_r\!\to\! 0$).
We verify these predictions empirically in \cref{sec:exp:subsec:exp3_model_validation}.

\subsection{PS Collapse}\label{sec:POPS:subsec:PS_collapse}
While PS is the dominant solution concept in the DD literature, it can be suboptimal relative to PO \cite[Prop.~2.1]{miller2021outside}.
In ONP, PS degenerates entirely.

\begin{theorem}[PS Collapse]\label{thm:PS_failure}
For ONP with $\lambda\!>\!0$, the unique PS point degenerates to $p_{PS} \!=\! \displaystyle \argmin_{p\in[p_{\min},p_{\max}]}\!\tfrac{\lambda}{2}\|p\|_2^2$, so $p_{PS}\!\to\! p_{\min}$ (or $\zeros$), ignoring congestion.
\end{theorem}
\noindent
PS collapse is a $\Pi$PO-specific pathology: in standard PO ($\Pi\!=\! I$), PS can be suboptimal \cite[Prop.~2.1]{miller2021outside} but does not degenerate to a trivial solution.
While in ONP, the failure is structural: PS decouples from the congestion objective regardless of congestion level.
\begin{itemize}[leftmargin=*]\setlength{\itemsep}{0pt}
\item \textbf{Conceptually:} a PS platform regards congestion as exogenous---it observes congestion but finds no utility in raising $p$, treating the distribution as fixed.
\item \textbf{Algorithmically:} PS-seeking methods
\cite{perdomo2020performative,drusvyatskiy2023stochastic,wang2025projected}
fail on ONP due to a structural mismatch (standard PS targets coupled bilinear losses $p^\top x$, whereas ONP is composite and non-bilinear) and the Gating Jacobian (PS ignores $\Pi_{[\zeros,x_{\max}]}$, so whenever $x\geq x_{\max}$ the update cannot ``see'' $\cD$).
\end{itemize}
Our approach fully accounts for the sensitivity of $\cD(p)$ to prices via the Gating Jacobian, yielding a Jacobian-aware $\Pi$PO-seeking algorithm.

\paragraph{When does PS equal $\mathit{\Pi}$PO}
Since $p_{PS}\!=\!p_{\min}$, coincidence requires $p^*_{\Pi\text{PO}}\!=\!p_{\min}$.
Inspecting \eqref{eq:KKT_true}, this forces the congestion term to vanish at $p\!=\!p_{\min}$.

\begin{proposition}[PS$=\!\Pi$PO Only in Degenerate Cases]\label{prop:PS_eq_PO}
For ONP with $\lambda\!>\!0$ and $\pi_{\min}\!>\!0$,
PS and $\mathit{\Pi}$PO coincide ($p_{PS}\!=\!p^*_{\Pi\text{PO}}\!=\!p_{\min}$)
only when:
(i) \textbf{trivial congestion} ($x_{\max}\!=\!\zeros$)
;
(ii) \textbf{no behavioral response} ($B\!=\!\zeros$); 
or 
(iii) \textbf{extreme regularization} ($\lambda\!\to\!\infty$).
In all other cases $p_{PS}\!\neq\! p^*_{\Pi\text{PO}}$.
\end{proposition}
\noindent
None of (i)--(iii) reflects a realistic network pricing scenario.
PS is not merely suboptimal but a \emph{structurally wrong} solution concept for ONP\footnote{The only continuous path from PS to $\Pi$PO is $\lambda\to\infty$, which destroys the congestion signal rather than making PS performatively aware.}: standard stability-seeking heuristics are inadequate for performative, capacity-constrained intelligent pricing systems, motivating $\Pi$\text{PO}-seeking for such environments.

\section{Trust-Region SQP with Dead Zone Escape}\label{sec:algo_design}
To solve the nonsmooth constrained SAA problem \eqref{prob:ONP-nox-SAA}, and enable an intelligent pricing system to reliably learn $\Pi$PO-optimal prices, we propose a TR-SQP algorithm \cite{nocedal2006numerical}.\footnote{Why TR-SQP:
\textbf{1.}\,Shrinking the TR radius $\Delta_t$ controls stepsizes near kinks, preventing oscillations across active-set boundaries in the partially gated landscape.
\textbf{2.}\,It handles stiff demand constraints $c_N(p) \!\leq\!\zeros$ via explicit linearization, avoiding ill-conditioning of penalty methods when CQ fails (Prop.~\ref{prop:fail_CQ}).
\textbf{3.}\,It leverages sparsity (\cref{sec:app:subsec:sparse}) for scalability.}
Standard SQP can fail due to \DZ; we design an escape mechanism to address this.

\paragraph{TR-SQP subproblem}
At each iteration $t$, we model $f_N$, $c_N$ around $p_t$ by a quadratic and linear model, respectively.
Let $H_t$ be a curvature matrix and $\Delta_t$ the trust-region radius; the step $d_t$ solves:\footnote{We use the infinity norm to keep \eqref{def:TR_model} a QP. The $\infty$-norm constraint naturally aligns with the box constraints $p \!\in\![p_{\min}, p_{\max}]$.}
\begin{equation}\label{def:TR_model}
\argmin_{\|d\|_\infty \!\leq\Delta_t}
m_t(d) ~\coloneqq~ \Big\langle \frac{1}{2}H_t d + \nabla_p f_N(p_t), ~d \Big\rangle
~~\st~~ c(p_t) \!\leq\!-\big\langle \nabla_p c(p_t), d \big\rangle.
\end{equation}
With $H_t\!=\!\nabla^2_p f_N$, the model $m_t$ is convex almost surely (Theorem~\ref{thm:fE_strong_cvx}).
See \cref{sec:app:subsec:sparse} for efficient computation of $(\nabla^2_p f_N)^{-1}$.

\paragraph{TR update and acceptance}
Assuming the iterate $p$ is outside the full \DZ, Algorithm~\ref{algo:TRSQP} uses the following merit functions.
\begin{itemize}[leftmargin=*]\setlength{\itemsep}{0pt}
\item \textbf{Inside active zone} ($\omega_{\mathrm{sum}} \!\geq\!\varepsilon_{\text{mach}}$):
steps are evaluated with the standard $\ell_1$ merit function
$\phi(p, \nu) \!=\! f_N(p) + \nu \|c_N(p)^+\|_1$ with $(x)^+ \!=\! \max\{\zeros,x\}$.
\item \textbf{Outside active zone ($p_r\in\mathrm{int}\DZ_r$)} ($\omega_{\mathrm{sum}} < \varepsilon_{\text{mach}}$):
steps use the \emph{augmented} merit function\footnote{Here we use $\mu=\IE[\zeta]$ instead of MC for smaller computation cost.} with $\mu_V \!>\!\lambda\|p_{\max}\|_2\|B^{-1}\|_\infty$:
\begin{equation}\label{def:aug_merit}
\Phi(p,\nu) \coloneqq \phi(p,\nu) + \mu_V V(p),
\quad
V(p) \coloneqq \bigl\| Bp+\mu - \Pi_{[0,x_{\max}]}(Bp+\mu) \bigr\|_1.
\end{equation}
\end{itemize}
Outside \DZ, $V(p)\!=\!0$ so $\Phi \!=\! \phi$.
The penalty parameter is updated as
$\nu_t \!=\! \min(\nu_{\max}, \max(\nu_{t-1}, \|\gamma\|_\infty \!+\! \delta))$.
Predicted and actual reductions are:
\[
\begin{array}{rl}
\mathrm{Pred}_t \hspace{-3mm}&= \nu_t \bigl(\mathrm{viol}_t - \|(c_N(p_t) + \langle \nabla_p c_N(p_t), d_t \rangle)^+\|_1\bigr) - \langle \tfrac{1}{2}H_t d_t + \nabla_p f_N(p_t), d_t \rangle, \\
\mathrm{Ared}_t \hspace{-3mm}&= \phi(p_t, \nu_t) - \phi(p_t + d_t, \nu_t),
\end{array}
\]
where $\mathrm{viol}_t\!=\!\|c_N(p_t)^+\|_1$.
Step acceptance follows the standard TR ratio $\rho_t \!=\!\mathrm{Ared}_t/\mathrm{Pred}_t$.\footnote{If $\rho_t\!>\!\eta_1$: accept ($p_{t+1} \!=\! p_t + d_t$); expand radius if $\rho_t\!>\!\eta_2$.
If $\rho_t\!\leq\!\eta_1$: reject ($p_{t+1}\!=\! p_t$), shrink radius.}
DZ escape steps are accepted unconditionally, as $\Phi$ decreases by construction (Lemma~\ref{lem:finite_dz_escape}).

\begin{algorithm}[H]
\small
\caption{TR-SQP with Dead Zone Escape}\label{algo:TRSQP}
Input $p_0 \!\in\! \IR^{|\cR|}, \lambda \!>\! 0,$ set param. $(\eta_1,\eta_2, \nu_0) \!=\! (0.25, 0.75, 1)$,
set $\mu_V \!>\! \lambda\|p_{\max}\|_2\|B^{-1}\|_\infty$ (offline),
set $\varepsilon_0 \!>\! \varepsilon_{\text{mach}}$ (fixed nudge size)
\\
\For{$t \!=\! 0, 1, 2, \ldots$}{
get $y^i \!=\! Bp_t+\zeta^i$ for all $i$
\hfill latent demand
\\
get $\omega^i = \II(\zeros \!<\! y^i \!<\! x_{\max})$ for all $i$
\hfill gating indicator
\\
    \uIf{$\omega_{\text{sum}} \!=\! \sum_{i=1}^N \|\omega^i\|_1  < \epsilon_{\text{mach}}$ (\textbf{Full DZ Detected})}{
    get $\overline{y} \!=\! \Avg[y^i]$ and 
    $\mathrm{res} \!=\! \overline{y} - \Pi_{[\zeros,x_{\max}]}(\overline{y})$ \hfill residual 
    \\
        \uIf{$\|\mathrm{res}\|_\infty \!=\! 0$ (measure-0 rare event)}{
        $p_{t+1} \!=\! p_t + B^{-1}\delta$,
        $
        \delta_r \!=\! \begin{cases}
        - \varepsilon_0 & \overline{y}_r \!=\! (x_{\max})_r
        \\[-3pt]
        +\varepsilon_0  & \overline{y}_r \!=\! 0
        \\[-3pt]
        0  & \text{else}
        \end{cases}
        $
        \hfill fixed nudge (rare)
        }
        \Else{
        $d_{\text{refl}} \!=\! -B^{-1} \mathrm{res}_t$
        \hfill reflection direction
        \\
        $d_t \!=\! \min(1, \Delta_t/\|d_{\text{refl}}\|_\infty) d_{\text{refl}}$
        \hfill scale with TR radius
        \\
        $p_{t+1} \!=\! p_t + d_t$
        \hfill escape by reflection
        }
        $\Delta_{t+1} \!=\! \Delta_t$        \hfill radius unchanged
    }
    \Else{
    get $J^i \!=\! \Diag(\omega^i)$, $x^i \!=\! \omega^i \odot y^i$
    \hfill Jacobian \& realized flow
    \\
    get $\nabla_p f_N$, $\nabla_p c_N$
    and $H_t \!=\! \nabla^2_p f_N$ (replace by $H_{GGN}$ if unstable)
    \\
        \uIf{Solving QP gives feasible solution}
        {
        $(d_t, \gamma)\!= $\textbf{QP Subproblem}
        \hfill standard SQP step
        }
        \Else{
        $\textstyle d_t \!=\! \argmin_{\|d\|_\infty \!\leq\!\Delta_t} \| c(p_t)^+ + \langle \nabla_p c, d \rangle \|_1, ~~\gamma \!=\! \zeros$
        \hfill Restoration
        }
        TR Update:
        Update $p_{t+1}, \gamma, \Delta_{t+1}$
    }
}
\end{algorithm}

\subsection{Algorithmic Behavior in Dead Zone Escape}\label{sec:algo_design:subsec:reflect}
Standard solvers fail in partial \DZ\ due to zero gradient (CQ failure, Prop.\ref{prop:fail_CQ}); our TR-SQP handles this 
via constraint restoration\footnote{The gated route $r$ is temporarily uncontrollable, while the restoration step remains well-defined on active rows.
Once $p$ escapes $\DZ_r$, it restores CQ, and re-enters the QP subproblem.
}.
The projection divides the landscape into active regions ($y \!\in\!(\zeros, x_{\max})$), strictly saturated interiors ($\mathrm{int}\DZ_r$), and measure-zero boundaries (kink $\partial\DZ$, where $\mathrm{res} \!=\! \zeros$).

By Corollary~\ref{cor:pIE_not_in_DZ}, $p^*_{\IE}\notin\mathrm{int}\DZ(f_N)$ 
w.h.p.; so it must lie in partial \DZ, a strictly active zone, or on a kink.
$f_N$ may have local minimizers inside $\mathrm{int}\DZ$; these are finite-sample artifacts 
that vanish as $N \!\to\! \infty$.
Algorithm~\ref{algo:TRSQP} evicts all iterates from $\mathrm{int}\DZ$ and targets Clarke-stationarity of $f_N$, which coincides with the global $\Pi$PO when $N \!\geq\! N_{\mathrm{op}}$.

\paragraph{Partial DZ and kinks ($p \!\in\! \DZ_r$ for some $r$)}
If some routes are gated ($J_{rr} \!=\! 0$ for some $r$), the remaining active routes still induce a nonzero gradient. 
If an iterate lands on a kink ($\partial\DZ$), TR radius shrinkage prevents oscillation, requiring no geometric override.

\paragraph{Full saturation ($p \!\in\!\mathrm{int}\DZ$, $J\!=\!\zeros$)}
When all routes are saturated, the gradient vanishes.
Two failure modes (of standard gradient-based solvers) arise:
\begin{itemize}[leftmargin=*]\setlength{\itemsep}{0pt}
\item \textbf{PS-Collapse ($c_N(p) \!\leq\!\zeros$):} Standard solvers follow $\nabla_p f_N(p) \!\overset{\eqref{eq:DZ_grad_collapse}}{=}\! \lambda p$, driving prices toward the trivial PS solution (Theorem~\ref{thm:PS_failure}).
\item \textbf{Saturation Trap ($c_N(p) \!>\! \zeros$):}
$J\!=\!\zeros \!\implies\! \nabla_p c_N \!=\! \zeros$ (Lemma~\ref{gradientHessian_derivation}), so restoration-based methods stagnate with CQ failure (Prop.~\ref{prop:fail_CQ}).
\end{itemize}
Algorithm~\ref{algo:TRSQP} detects full saturation via $\omega_{\mathrm{sum}}\!<\! \varepsilon_{\text{mach}}$ and triggers a \textbf{reflection step} (or, in the measure-zero boundary case, a \textbf{nudge step}).

\paragraph{Merit function for DZ steps}
To give DZ escape steps a Lyapunov certificate, we verify that the reflection $d_{\mathrm{refl}} \!=\! -B^{-1}\mathrm{res}$ is a descent direction for $\Phi$:
\[
\begin{array}{rcl}
\langle \nabla_p V, d_{\mathrm{refl}}\rangle
&=&
\mathrm{sgn}(\mathrm{res})^\top B(-B^{-1}\mathrm{res})
= -\|\mathrm{res}\|_1 < 0,
\\
\langle \nabla_p \Phi, d_{\mathrm{refl}}\rangle
&\leq&
\lambda\|p_{\max}\|_2\|B^{-1}\|_\infty\|\mathrm{res}\|_1 - \mu_V\|\mathrm{res}\|_1
\\
&=&
-(\mu_V - \lambda\|p_{\max}\|_2\|B^{-1}\|_\infty)\|\mathrm{res}\|_1 < 0.
\end{array}
\]
So every reflection step strictly decreases $\Phi$.
As $\Phi$ is bounded below (compact feasible set, $V\!\geq\!0$), DZ re-entry occurs finitely many times across the entire run\footnote{Each DZ step decreases $\Phi$ by at least $c_\Phi \!>\!0$.
Since $\Phi \!\geq\!\Phi_{\min}$, total DZ steps are at most $(\Phi_0 - \Phi_{\min})/c_\Phi$.
See Lemma~\ref{lem:finite_dz_escape}.}.
The nudge step also decreases $\Phi$ and terminates in a single step (Lemma~\ref{lem:finite_dz_escape}).

\paragraph{GGN Hessian}
In regions where $J^i\!\approx\!\zeros$ across most samples, $\nabla^2_p f_N\!\approx\!\lambda I$ causes ill-conditioning.
We use the approximation $H_{\text{GGN}}\!\coloneqq\! B^\top Q B + \lambda I$.\footnote{$H_{\text{GGN}}$ acts as a Levenberg--Marquardt regularizer, computed once offline, and is asymptotically exact as $\pi \!\to\! 1$: $H_{\text{GGN}} - \nabla^2_p f_N \!=\! B^\top \Avg[Q - J^i Q J^i] B$ with $\IE[(Q - J^i Q J^i)_{rs}] \!=\! (1 - \pi_r \pi_s) Q_{rs}$.}

\subsection{Global Convergence}\label{sec:algodesign:subsec:convergence}
Nonsmoothness in ONP invalidates standard SQP convergence theory; we prove Algorithm~\ref{algo:TRSQP} globally converges to a Clarke-stationary point \cite{clarke1975generalized}.

\begin{lemma}[Finite \DZ\ Escape]\label{lem:finite_dz_escape}
Let $\mu_V\!>\!\lambda\|p_{\max}\|_2\|B^{-1}\|_\infty$.
When Algorithm~\ref{algo:TRSQP} enters \DZ, one of two cases applies almost surely:
\\
\textbf{(i)} $\|\mathrm{res}_t\|_\infty\!>\!0$ (a.s. case):
Reflection steps yield strict contraction
$\|\mathrm{res}_{t+1}\|_1\!\leq\! (1\!-\!\alpha_t)\|\mathrm{res}_t\|_1$, strict decrease $\Phi(p_{t+1},\nu)\!\leq\!\Phi(p_t,\nu) \!-\! c_\Phi$, where constant $c_\Phi \!=\! (\mu_V \!-\! \lambda\|p_{\max}\|_2\|B^{-1}\|_\infty)\alpha_t\|\mathrm{res}_t\|_1 \!>\!0$,
and escape in at most $T_{\DZ} \!\leq\!\lceil\|\mathrm{res}_0\|_1/c_{\mathrm{dec}}\rceil$ steps where $c_{\mathrm{dec}} \!=\! \Delta_{\mathrm{entry}}/\|B^{-1}\|_\infty$.
\\
\textbf{(ii)} $\|\mathrm{res}_t\|_\infty\!=\!0$ (measure-zero case):
A single nudge step of size $\varepsilon_0\!>\!\varepsilon_{\text{mach}}$ moves $\bar{y}_r$ strictly off the boundary, after which the algorithm either exits \DZ\ or enters case (i).
\\
In both cases, escape is finite almost surely, and each step strictly decreases $\Phi$.
\end{lemma}
\noindent
At the post-escape iterate $p'$, $\pi_{\min}(p') \!>\!0$, restoring the gradient signal and making $N_{\DZ}$ finite (Remark~\ref{rem:N1_algo}).

\begin{theorem}[Global Convergence to Clarke-Stationarity]\label{thm:global_convergence}
Assume $\{H_t\}$ is uniformly bounded and the constant 
$\mu_V \!>\! \lambda\|p_{\max}\|_2\|B^{-1}\|_\infty$.
Then we have every accumulation point of $\{p_t\}$ is Clarke-stationary for $\phi$.
\end{theorem}

\paragraph{Math--algorithm correspondence}

As stated in Remark~\ref{rem:N1_algo}, $N_{\DZ}, N_{SAA}, N_{\mathrm{op}}$ diverge as $\pi_{\min} \!\to\! 0$ in \DZ.
The DZ escape (Lemma~\ref{lem:finite_dz_escape}) eliminates this cost by actively moving $p$ to where $\pi_{\min} \!>\!0$, after which all bounds become finite.
Theorem~\ref{thm:global_convergence} ensures post-escape iterates converge to Clarke-stationarity, which with $N \!\geq\!N_{\mathrm{op}}$ coincides with the unique global $\Pi$PO (Theorem~\ref{thm:fE_strong_cvx}), with the following approximation guarantee:

\begin{corollary}[End-to-End Guarantee]\label{cor:end_to_end}
Let $\lambda \!>\!0$, $\delta \!\in\!(0,1)$, $\epsilon_{opt} \!>\!0$, and $N \!\geq\!N_{\text{op}}$.
With probability at least $1\!-\!\delta$, Algorithm~\ref{algo:TRSQP} outputs $p_N^*$ satisfying $\|p_N^* - p^*_{\IE}\| \!\leq\!\tfrac{\epsilon+\epsilon_{\text{opt}}}{\alpha}$,
where $\alpha \!=\! \lambda + \pi_{\min}^2\lambda_{\min}(Q)\sigma_{\min}^2(B)$.
\end{corollary}

\begin{remark}[Total Solution Complexity]
For TR iteration count $T \!=\! \cO(\epsilon_{\mathrm{opt}}^{-2})$
\cite{nocedal2006numerical} and sample complexity $N_{\mathrm{op}}$, 
the total cost to reach $\epsilon_{\text{opt}}$-stationarity is
\[
T_{\mathrm{total}}
= \underbrace{\cO\bigl(\epsilon_{\mathrm{opt}}^{-2} \bigr)}_{\text{iteration cost}}
\times
\underbrace{\Tilde{\cO}\bigl(\pi_{\min}^{-4}|\cR|\ln|\cR|/\delta\bigr)}_{N_{\mathrm{op}} \text{ samples}}
~~+ \hspace{-0.1cm}
\underbrace{\cO(T_{\DZ})}_{\text{finite DZ-escape  (Lemma~\ref{lem:finite_dz_escape})}}.
\]
This bound characterizes the total data and computation needed to learn near-optimal prices under performative, capacity-constrained feedback.
\end{remark}

\paragraph{Summary: 2-World}
Fig.\,\ref{fig:two_worlds} summarizes the relationship between $\ftrue$ and $f_N$.

\begin{figure}[h!]
\centering\includegraphics[width=0.99\linewidth]{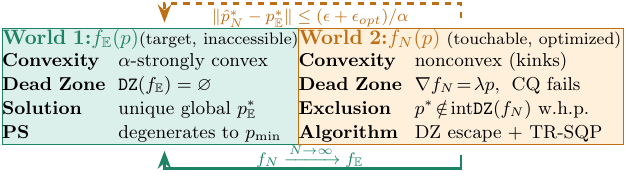}
\caption{The algorithm operates on $f_N$ (nonconvex, contains Dead Zones) but targets the unique global $\Pi$PO $p^*$ of $\ftrue$ ($\alpha$-strongly convex, $\DZ(\ftrue)\!=\!\varnothing$).
Corollary~\ref{cor:pIE_not_in_DZ} ensures $p^*\notin\mathrm{int}\DZ(f_N)$ w.h.p.;
the DZ escape (Lemma~\ref{lem:finite_dz_escape}) handles transient iterate traps;
uniform concentration $f_N\!\to\!\ftrue$ (Prop.~\ref{prop:concentration_combined}) bridges the two landscapes;
Corollary~\ref{cor:end_to_end} guarantees $\|p_N^* - p^*_{\IE}\| \!\leq\!\tfrac{\epsilon+\epsilon_{\text{opt}}}
{\alpha}$ for $N\geq N_{\mathrm{op}}$.}
\label{fig:two_worlds}
\end{figure}
 
\section{Experimental Validation}\label{sec:exp}
We validate our theoretical contributions, algorithm robustness, and scalability, demonstrating that our framework reliably learns near-optimal prices under performative, capacity-constrained feedback.
\cref{sec:exp:subsec:exp1_tqsrp} shows that TR-SQP finds the global $\Pi$PO that strictly dominates PS, and that 1st-order methods fail due to CQ violations.
\cref{sec:exp:subsec:exp2_robustness} validates DZ escape and convergence under pathological initializations.
\cref{sec:exp:subsec:exp3_model_validation} validates the ONP model and demonstrates that exact Jacobian treatment is strictly necessary.

\paragraph{Braess network and Social Optimum}
We use the Braess graph~\cite{braess1968paradoxon} (Fig.\,\ref{fig:Braess_graph}), a fundamental example for paradoxical congestion effects, as baseline to verify ONP optimality conditions.
With demand $l\!=\!1$, the Social Optimum (SO) splits flow equally between $r_1, r_3$ ($x_{r_1}^{SO}\!=\!x_{r_3}^{SO}\!=\!0.5$, avoiding the zero-cost bridge $r_2$), giving cost $\ell^{SO} \!=\!1.5$ versus Nash Equilibrium cost $\ell^{NE} \!=\! 2$.

\paragraph{Baselines and setup}
Unlike the deterministic Braess model, ONP's objective is analytically intractable.\footnote{Due to stochasticity and $\Pi_{[0, x_{\max}]}$; $\IE[\Pi(y)] \neq \Pi(\IE[y])$.}
We evaluate TR-SQP on \eqref{prob:ONP-nox-SAA} against:
\textbf{Brute Force} (numerical ground truth);
\textbf{Repeated Risk (RR) Minimization}~\cite{perdomo2020performative} (PS-seeking)\footnote{At iteration $t$: $p_t \!=\! \text{argmin}_{p} \text{cost}_{\text{tot}}(p, \cD(p_{t-1}))$.};
\textbf{1st-order methods} (GD/SGD); and \textbf{APO} (Def.~\ref{def:approxPO}).
We configure a highly congested regime to test capacity gating and DZ escape.
See Appendix~\ref{app:exp_setup} for details.

\subsection{Exp1: Failure of PS and 1st-order methods}\label{sec:exp:subsec:exp1_tqsrp}

\begin{figure}[ht]
\centering
\setlength{\tabcolsep}{2pt}
\renewcommand{\arraystretch}{0.85}
\raisebox{-0.45\height}{
\includegraphics[width=0.175\linewidth]{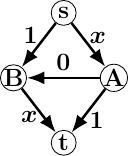}
}
\begin{tabular}{l|c|cccc|c}
& \textbf{SO} & \textbf{BF} & \textbf{RR(PS)} & \textbf{SGD} & \textbf{GD} & \textbf{TR-SQP} \\ 
\hline
Price $p_1$($\!=\!p_3$) & 0.5 & \textbf{4.749} & 0.002 & 5.058 & 5.069 & \textbf{4.751} \\
Price $p_2$ & 1.0 & 5.092 & 0.002 & 0.000 & 5.081 & \textbf{5.048} \\
Flow $x_1$($\!=\!x_3$) & 0.5 & \textbf{0.500} & 1.000 & 0.000 & 0.000 & \textbf{0.500} \\
Flow $x_2$ & 0.0 & \textbf{0.000} & 1.000 & 1.000 & 0.000 & \textbf{0.001} \\
\textbf{Objective} & - & 1.541 & 10.00 & $\infty$ & $\infty$ & \textbf{1.540} \\ 
\textbf{Time(s)} & - & 5.392 & 0.008 & - & - & \textbf{0.016} \\
\textbf{Gap\%} & - & - & 549.1 & $\infty$ & $\infty$ & \textbf{-0.06} 
\end{tabular}
\centering
\includegraphics[width=\linewidth]{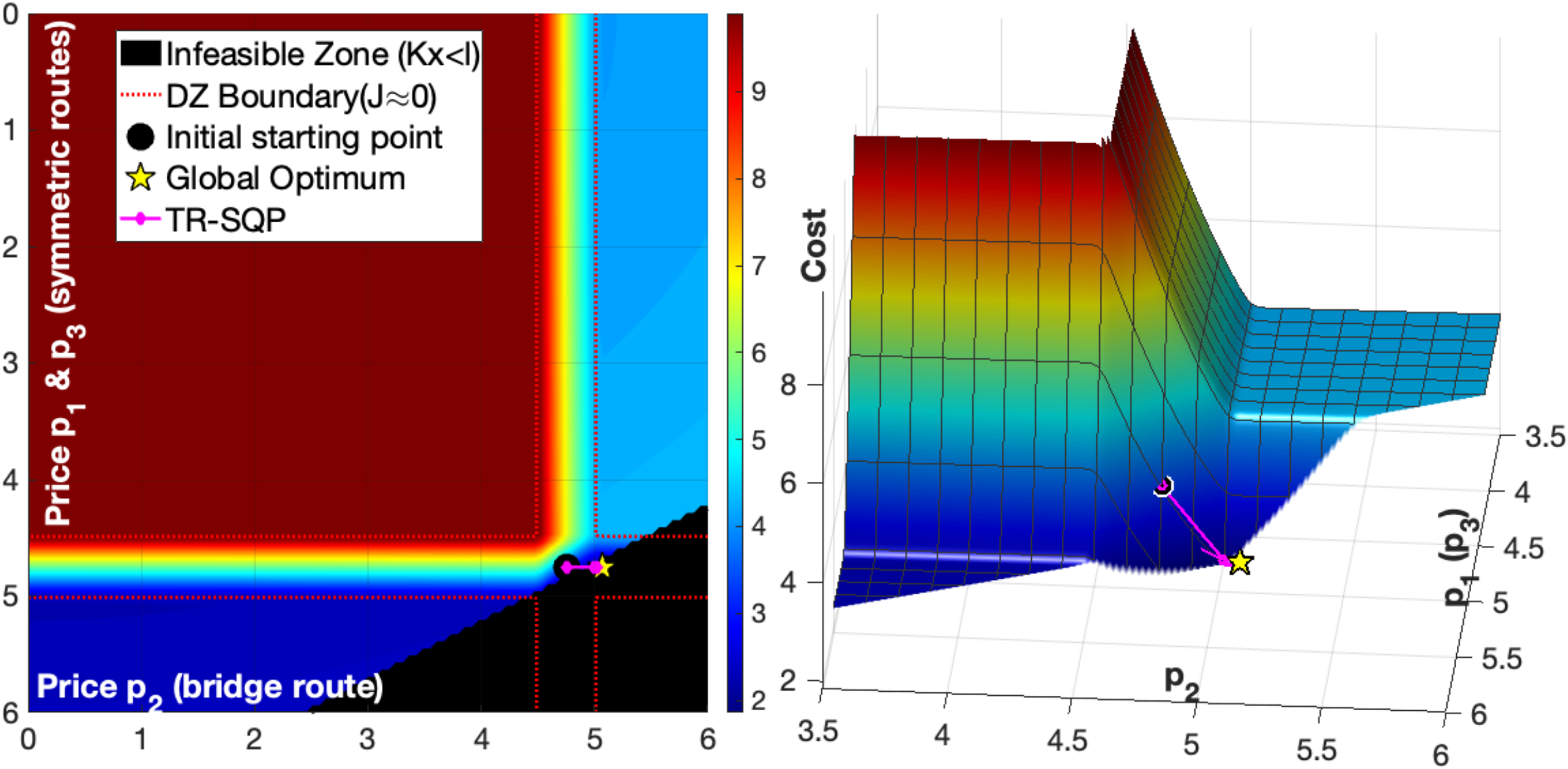}
\caption{Top: Results on the Braess network (`$\infty$' denotes constraint violation).
$\text{Gap\%} \!=\! (f_N(p) \!-\! f_N(p_{\text{BF}}))/ f_N(p_{\text{BF}})\times\!100\%$.
Bottom left: top-view of the objective with infeasible zone (black) and dead zone (red-dotted corners).
Bottom right: 3D surface zoomed near the global optimum.
Solid dot: starting point; purple line: TR-SQP trajectory converging to the optimum (star).
The contour plot reveals a narrow feasible corridor for algorithmic progress.}
\label{fig:Braess_graph}
\end{figure}

\paragraph{Failure of 1st-order methods and PS}
GD/SGD diverge, validating Prop.~\ref{prop:fail_CQ}: CQ failure in \DZ\ makes gradient methods structurally inapplicable.
RR converges but to near-zero prices ($p \approx 0.002$), distributing flow across all routes, including the inefficient bridge $r_2$, yielding a cost $10$ versus the optimum $\approx\! 1.5$.
This confirms Theorem~\ref{thm:PS_failure}: PS is blind to the full DD map.

\paragraph{TR-SQP performance}
TR-SQP correctly prices the inefficient bridge high ($p_2\!\approx\! 5$) and efficient routes lower ($p_1, p_3\!\approx\! 4.7$), steering flow to the optimal split $x\!\approx\! [0.5, 0, 0.5]^\top$.
It achieves a lower objective ($1.540$ in $0.016$s) than Brute Force ($1.541$ in $5.392$s), as Brute Force is grid-limited while TR-SQP performs continuous optimization, confirming Theorem~\ref{thm:fE_strong_cvx}.
Fig.\,\ref{fig:Braess_graph} shows the global $\Pi$PO lies in a narrow canyon adjacent to the infeasible zone; the TR mechanism prevents overshooting into \DZ\ or infeasibility.

\subsection{Exp2: DZ Escape and Boundary Convergence}\label{sec:exp:subsec:exp2_robustness}
We validate global convergence (Theorem~\ref{thm:global_convergence}) and the DZ escape mechanism (Lemma~\ref{lem:finite_dz_escape}) under three representative pathological initializations: full infeasible zone, Full DZ, and Partial DZ.
With the same setup as in Exp1 (\cref{sec:exp:subsec:exp1_tqsrp}), the DZ boundaries partition price space: $p\!<\!4.5$ triggers upper \DZ, $p\!>\!5$ triggers lower \DZ.
In the top two rows in Fig.\,\ref{fig:exp2traj}, all initializations converge to the unique global optimum.
We see that for the cyan curve, its 2nd iteration hits the kink, which shows that our TR-SQP is robust to the edge case.
The second row of Fig.\,\ref{fig:exp2traj} shows the result on a more complicated elasticity matrix $B$.

\paragraph{Convergence at kinks}
To test convergence when the optimum lies on a kink, 
we design two pathological cost structures\footnote{(1) prohibitive outer‑edge costs $c_2\!=\!c_4\!=\!100$, forcing the optimal flows $x_1\!=\!x_3\!=\!0$ (Dead Sides); (2) a prohibitive bridge cost $c_3\!=\!100$, forcing $x_2\!=\!0$ (Dead Bridge). For each case, TR-SQP parameters are tuned, while the problem's $\lambda$ remains unchanged.}.
The last row in Fig.\,\ref{fig:exp2traj} shows that the TR-SQP converges to the same kink-optimum: TR radius shrinks near the kink, providing a damping effect that stabilizes the steps near it.

\begin{figure}[h!]
\centering
\includegraphics[width=\linewidth]{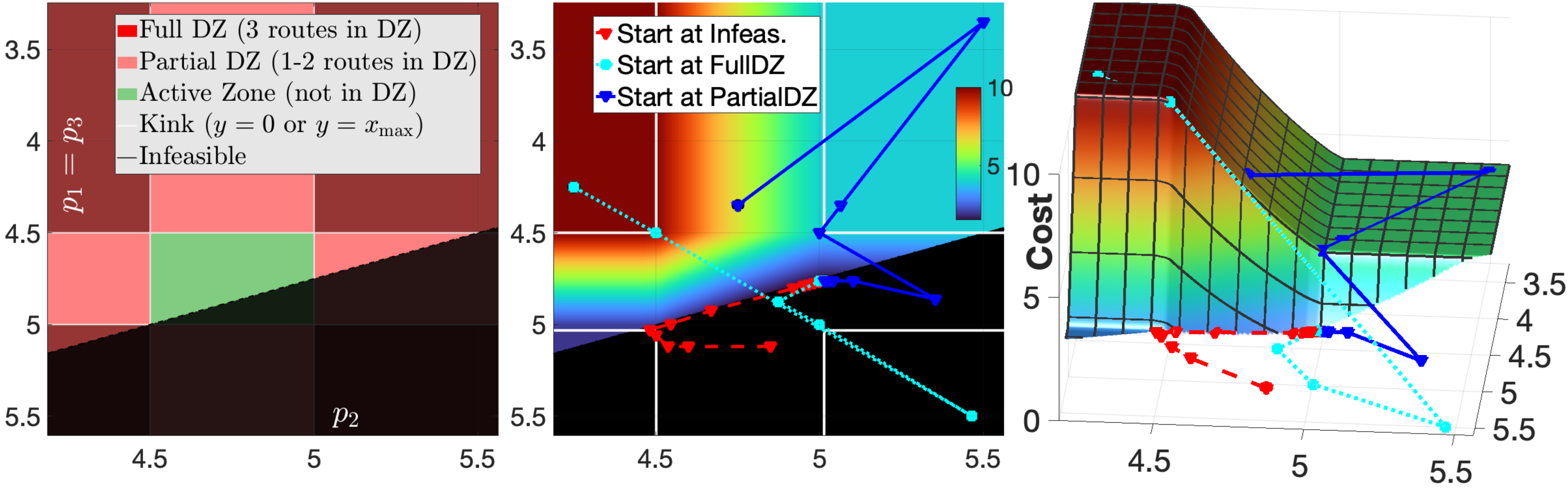}
\includegraphics[width=\linewidth]{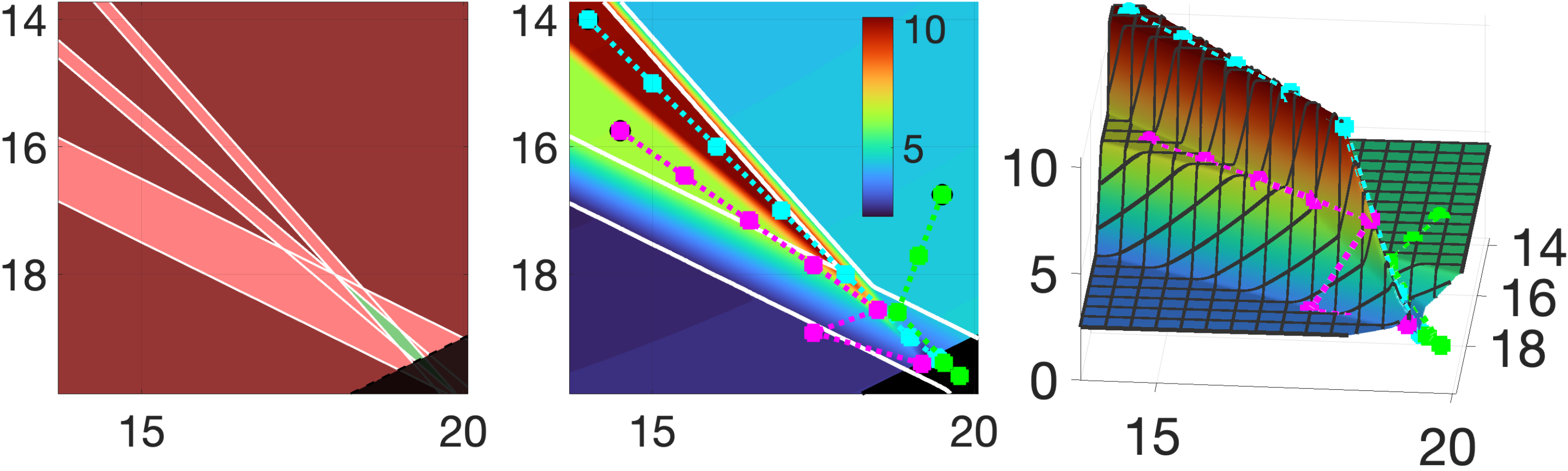}
\includegraphics[width=0.95\textwidth]{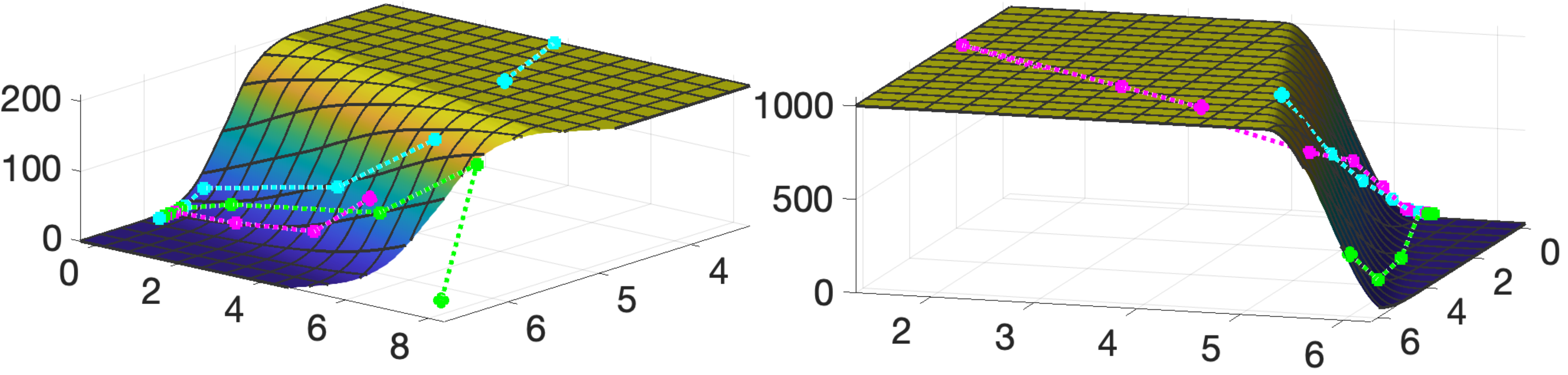}
\caption{\textbf{Top:} TR-SQP trajectories of three initializations.
\textbf{Mid:}  experiment with a more complicated $B$ encoding substitution effect.
\textbf{Bottom:} TR-SQP trajectories across three initializations where the $\Pi$PO lies on $\partial\DZ$, all converge to the Clarke-stationary $\Pi$PO.
}
\label{fig:exp2traj}
\end{figure}

\subsection{Exp3: ONP Model Validation}\label{sec:exp:subsec:exp3_model_validation}
We validate the performative flow steering (\cref{sec:algoframeworkSAA:subsec:stearing}) and quantify the necessity of exact Jacobian treatment.

\paragraph{Flow steering}
Fig.\,\ref{fig:DDDsterring} visualizes how TR-SQP steers $\cD(p_t)$ toward the $\Pi$PO.
We initialize with $p_0\!=\!B^{-1}(0.9\cdot\ones\!-\!\mu)$ (deeply saturated regime) and use a small TR radius ($\Delta_{\mathrm{init}} \!=\!0.001$, $\Delta_{\max}\!=\!0.1$) to make distributional drift visible.

\begin{figure}[h!]
\centering
\includegraphics[width=\linewidth]{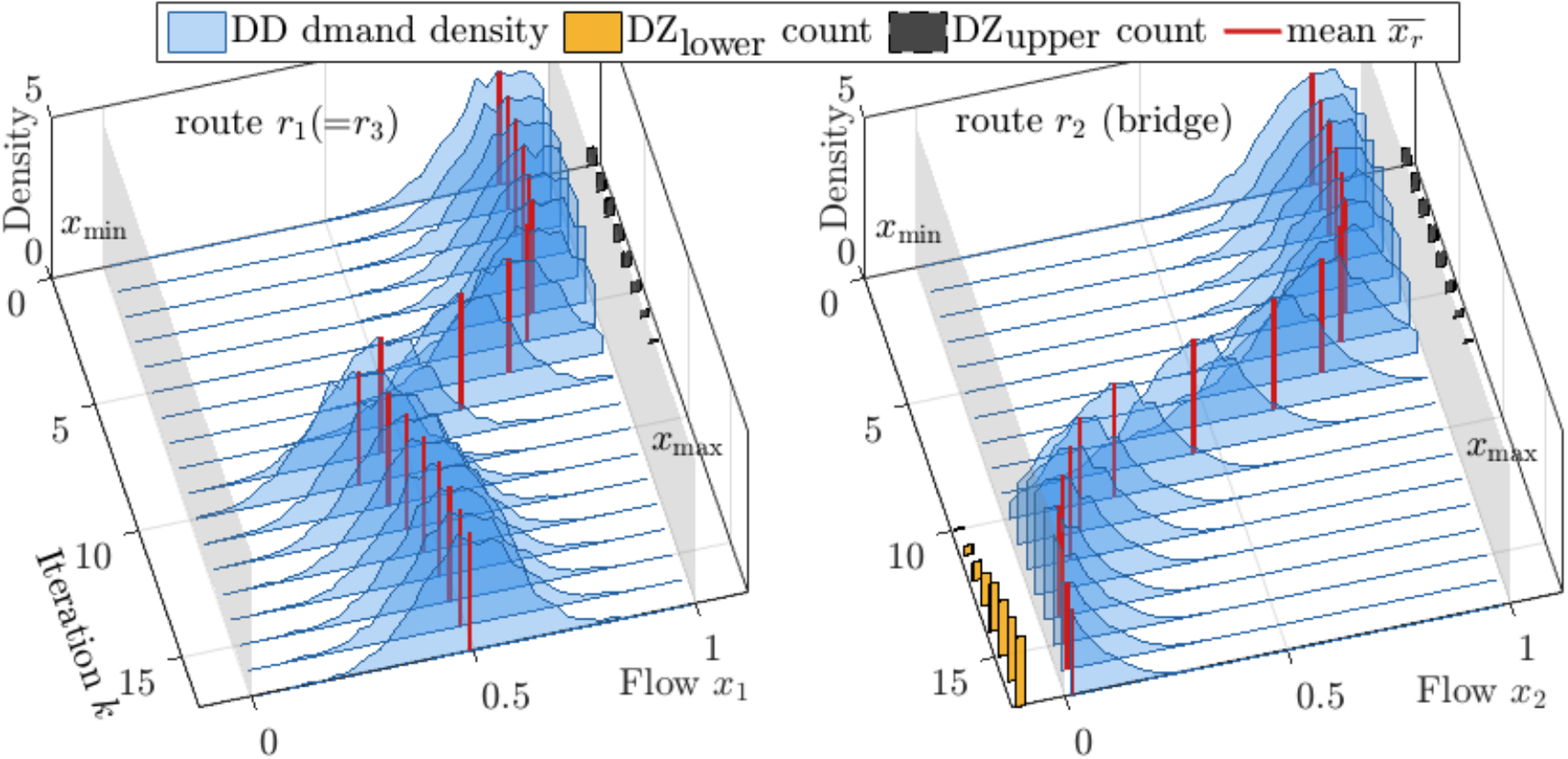}
\caption{DDD map $x_r \!\sim\! \cD(p_t)$ shifting as TR-SQP updates $p_t$ across iterations $k\!=\![0,17]$.
\textbf{Left} ($r_1 \!=\! r_3$): mean flow moves from $\bar{x}_1 \!\approx\! 0.9$ toward the SO at $\bar{x}_1^* \!\approx\! 0.5$; upper DZ mass shrinks as the distribution exits saturation.
\textbf{Right} ($r_2$): TR-SQP raises the price, pushing the distribution toward $x_{\min} \!=\! \zeros$; lower DZ mass grows, reflecting $J_{22} \!\to\! 0$.
Each price update shifts $\cD(p_t)$, which alters $\nabla f_N(p_t)$ for the next step.}
\label{fig:DDDsterring}
\end{figure}

\paragraph{APO-$\Pi$PO gap}
We quantify $\texttt{Gap}\!=\!|f_N(p_{APO}) \!-\! f_N(p_{\Pi\text{PO}})|$ by sweeping $x_{\max}$ and $\lambda$ on the Braess network under high latent demand (see \cref{app:exp_setup_exp3}).
Fig.\,\ref{fig:Exp3} confirms the regime-dependent predictions:
descending $x_{\max}$ widens \texttt{Gap} at high $\lambda$ (Prop.~\ref{prop:reg_amp_gap});
ascending $x_{\max}$ widens \texttt{Gap} at low $\lambda$ (Prop.~\ref{prop:reg_amp_gap});
at the equilibrium $\lambda^*$, \texttt{Gap} remains near zero across $x_{\max}$ (Prop.~\ref{thm:global_gap_amp}).
Crucially, APO becomes \emph{infeasible} in congested regimes (gray region in Fig.\,\ref{fig:Exp3}): it fails precisely when the exact Jacobian treatment matters most.

\begin{figure}[h!]
\centering
\includegraphics[width=\linewidth]{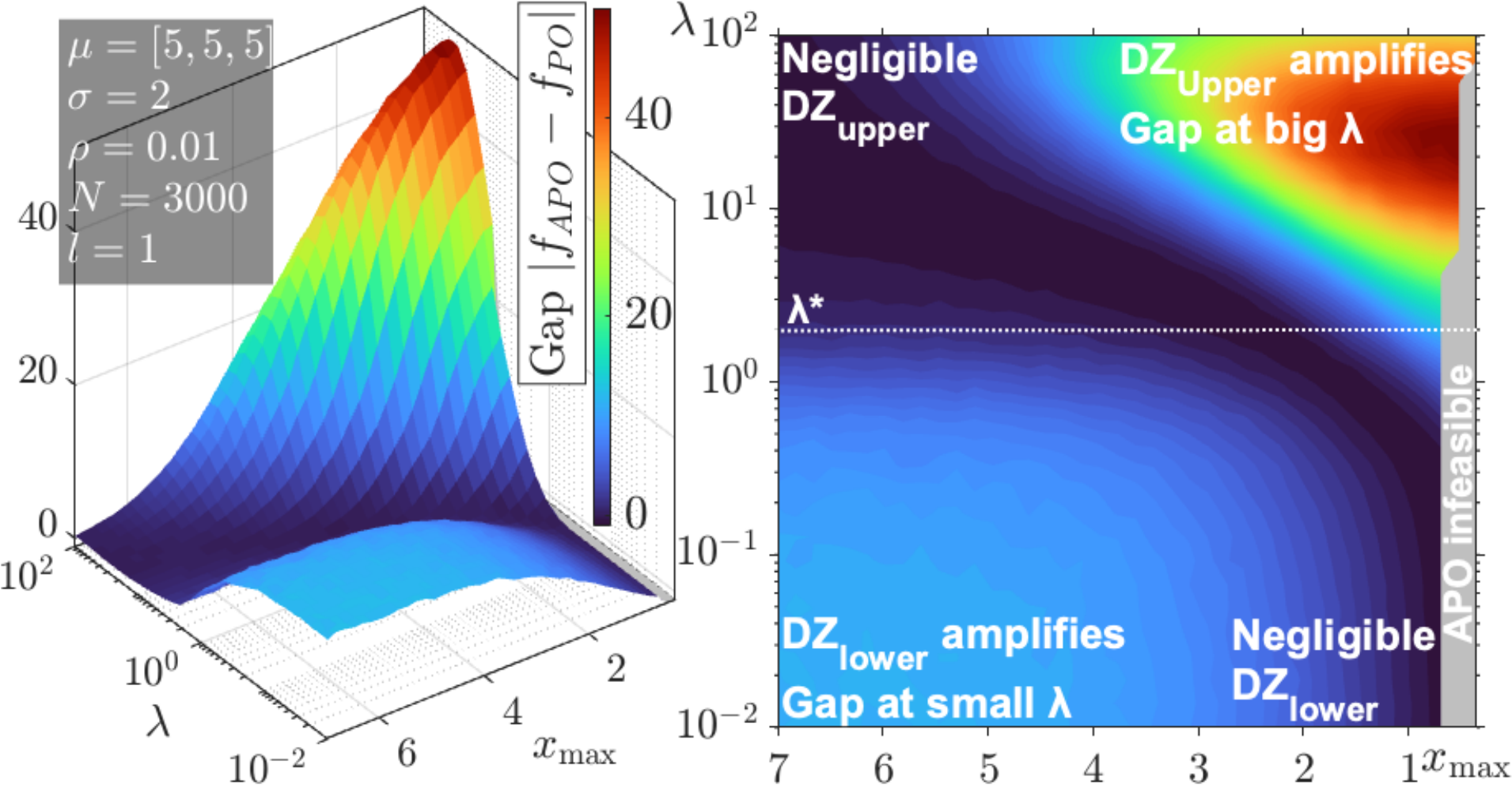}
\caption{Left: \texttt{Gap} surface across $x_{\max}$ and $\lambda$.
Right: 2D contour. Gray marks APO infeasibility due to constraint violations.}
\label{fig:Exp3}
\end{figure}

\begin{remark}[Equivalent effects]
Similar asymmetric DZ effects arise from varying demand $l$: raising $l$ pushes latent demand into upper \DZ\ (high-$\lambda$ gap inflation); lowering $l$ relaxes toward lower \DZ.
At $\lambda^*$, the gap can still be widened by globally reducing $\pi_r$ (e.g., increasing $\sigma$), consistent with Prop.~\ref{thm:global_gap_amp}.
\end{remark}

Together, these experiments confirm that correctly modeling the gating Jacobian is a practical necessity for any intelligent pricing systems operating near capacity boundaries.
\section{Scalability and Generalization}\label{sec:scalability}
We move to other networks to validate scalability on real-world topologies, demonstrating that our framework scales to the large networks where AI pricing systems operate in practice.

\subsection{Graph-Aware Scalability}\label{sec:exp:subsec:smallworld}
Real-world graphs are large and thus the sample complexity bounds in \cref{sec:DZ_SampleComplexity} can be huge.
However, real-world networks exhibit planarity and small-world structure that induce \textit{double sparsity}, substantially tightening these bounds.
The small-world effect \cite{watts1998collective} gives average path length $\overline{L} \!\approx\! \cO(\ln |\cV|)$, so the assignment matrix density decays:
$\nnz(A)/(|\cR| |\cE|)
\!\approx\! \tfrac{\overline{L}}{|\cE|}
\propto \tfrac{\ln |\cV|}{|\cV|} \!\to\! 0$ as $|\cV| \!\to\! \infty$.
This reduces per-iteration gradient and Hessian cost (\cref{sec:app:subsec:sparse}).

\paragraph{Graph-aware activation probability}
Our DDD model accumulates uncertainty along routes: $\Sigma\!=\! \sigma^2(I + \rho AA^\top)$.
The activation probability decays with route length as $\pi_r \!\propto\! 1/\sqrt{(AA^\top)_{rr}}$,\footnote{Probability mass scales inversely with standard deviation $\pi_r \!\propto\! \tfrac{1}{\sigma_r}$ and $\sigma_r \!\propto\! \sqrt{(AA^\top)_{rr}}$.} encoding a ``weakest link'' effect: a route is dead-zoned if \textit{any} edge saturates, so $\pi_r$ decays as the route lengthens.

\paragraph{Graph-Aware Sample Complexity}
If the graph is planar with maximum route length $L_{\max} \!\propto\! |\cR|^\gamma$, the operational sample complexity can be dropped to $N_{\text{graph}} \!=\! \tilde{\cO}(|\cR|^2\ln |\cR|)$\footnote{As $\sqrt{(AA^\top)_{rr}} \!=\! \sqrt{L_r}$.
With $L_{\max}\!\propto\!|\cR|^\gamma$: $\pi_{\min} \!\propto\! |\cR|^{-\gamma/2}$.
Substituting into $N_{\text{op}}$, assuming $\gamma\!=\! 2$ for planar graph.}, so finding the $\Pi$PO is tractable for large networks.

\subsection{Exp4: Simulation on the real-world GÉANT Network}\label{sec:exp:subsec:geant}
We validate our framework on the \textbf{GÉANT Network}~\cite{sndlib}, a pan‑European IP backbone.

\begin{figure}[ht!]
\raisebox{-0.5\height}{
\includegraphics[width=0.205\linewidth]{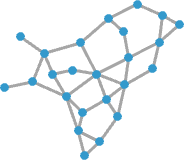}
} 
\begin{tabular}{l c c}
Metric & $\Pi$PO (ours) & APO \\
\hline 
Final cost & $1.57 \times\! 10^6$ & $8.41 \times\! 10^6$ \\
constraint violation & $3.64 \times\! 10^{-3}$ & $1.68 \!\times\! 10^{-2}$ \\
DZ fraction at solution & 64.1\% & 90.9\% \\
Mean $\bar\pi_r$ / Min $\pi_{\min}$ & 0.359 / $\approx\! 0$ & --- \\
$\mathrm{nnz}(A) / (|\cR||\cE|)$ & 4.78\% & --- 
\end{tabular}
\caption{
Results on GÉANT network (24 nodes, 74 directed edges) with $|\cR|\!=\!1510$ routes ($\leq\! 3$ paths per commodity, consistent with A1), mean path length $\bar{L}=3.53$, 538 commodities, $N\!=\!5000$, volume-to-capacity ratio $0.85$.
Constraint violation is defined as $\max_k (l_k - [K\bar x]_k)_+$.}
\label{fig:GEANT}
\end{figure}

\paragraph{Elasticity matrix generation}
We set $B \!=\! \tfrac{1}{2}(K^\top K) \odot \texttt{rand} \!-\! 2I$, where $\texttt{rand}\!\in\! [0,1]^{|\cR| \times |\cR|}$ has i.i.d.\ entries from the uniform distribution $\cU(0,1)$.
The term $K^\top K$ encodes route overlap (shared edges) and induces commodity-induced sparsity, while $\texttt{rand}$ adds heterogeneous  cross-elasticities for economic realism.

We run three trials with different random seeds, obtaining consistent results (within $\pm 3\%$); we report one representative run below.

\paragraph{Optimality}
$\Pi$PO achieves cost $1.57\!\times\! 10^6$ with negligible constraint violation ($3 \times\! 10^{-3}$).
APO incurs $434\%$ higher cost with 10 times larger violation, showing the Jacobian is necessary at scale.

\paragraph{Dead Zone}
At the initial starting point, only 7.1\% of routes lie in \DZ\ ($J_{rr}\!=\!0$), with $\overline{\pi}_r=0.9$ ($\pi_{\min}=0.76$).
One might be tempted to dismiss DZ handling as unnecessary given this small starting fraction.
However, after running the solver, the $\Pi$PO solution places $64.1\%$ of routes lie in \DZ, revealing substantial saturation.
This stark increase confirms that our TR-SQP algorithm correctly incorporate the nonsmooth Jacobian and DZ effects, far from being negligible, dominate the final solution.
Here $\pi_{\min} \!\approx\! 0$ indicates near-permanently saturated routes, confirming DZ escape (Lemma~\ref{lem:finite_dz_escape}) remains necessary.

\paragraph{Sparsity}
$\tfrac{\nnz(A)}{|\cR||\cE|} \!=\! 4.78\%$ with $\overline{L} \!=\! 3.53$, consistent with the small-world prediction $\ln_2 24 \!\approx\! 4.5$ (\cref{sec:exp:subsec:smallworld}).
This double sparsity reduces per-iteration cost.
\subsection{Exp5: Social Media Attention Allocation}\label{sec:exp:subsec:Twitch}
ONP also applies to flow control in digital ecosystems beyond transportation and IP networks.
In recommendation systems on YouTube/TikTok, the platform acts as a mechanism designer, shaping user behavior and allocating audience attention across creators while avoiding viewer burnout. 
ONP captures this: $p$ is an incentive signal (e.g., ad load), steering content generation and distribution (flow $x$) across demographic or interest-based ``routes''.

In a social network topology, $Q$ penalizes information redundancy and echo-chamber effects \cite{cinelli2021echo}: when excessive information flows through highly overlapping social circles, its marginal value decreases due to congestion.  
This enables the framework to dynamically mitigate attention saturation.
We validate this on the Twitch \texttt{ENGB} social network, where nodes are streamers and flows are followers.
Here, ONP allocates viewer attention across channels to minimize aggregate fatigue under content consumption demands.

\paragraph{Efficiency and Scalability}
We validate TR-SQP's computational efficiency by comparing the wall-clock runtime of dense vs sparse implementations (details in \cref{sec:app:subsec:exp5}).
Fig.\,\ref{fig:exp5:all} shows that the sparse solver scales better with $|\cR|$, achieving up to $4\times$ speedup at $|\cR|\!\approx\! 10^{4}$. 
The advantage diminishes at $|\cR|\!\approx\! 15000$ as the shared QP sub-solver becomes the dominant cost.
This confirms \cref{sec:exp:subsec:smallworld} that exploiting sparsity 
$(\tfrac{\nnz}{|\cR||\cE|}\!\approx\! 0)$ makes ONP tractable for large-scale deployments.
Furthermore, $\Pi$PO solvers reduce cost by $60\%$ versus APO at large scales. 
These results confirm that $\Pi$PO-targeting is theoretically sound and practically deployable across diverse large-scale networked environments.
\begin{figure}[h!]
\centering
\includegraphics[width=0.52\linewidth]{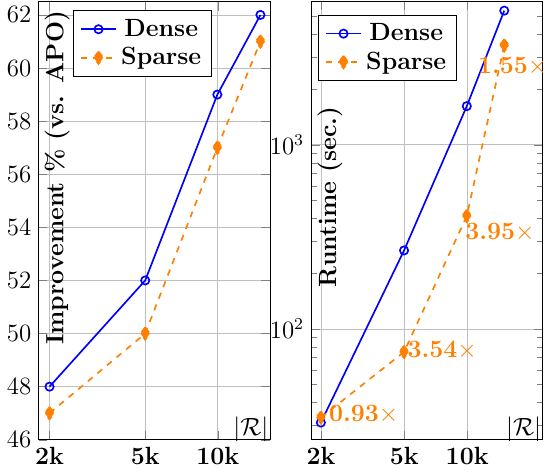}
\includegraphics[width=0.44\linewidth]{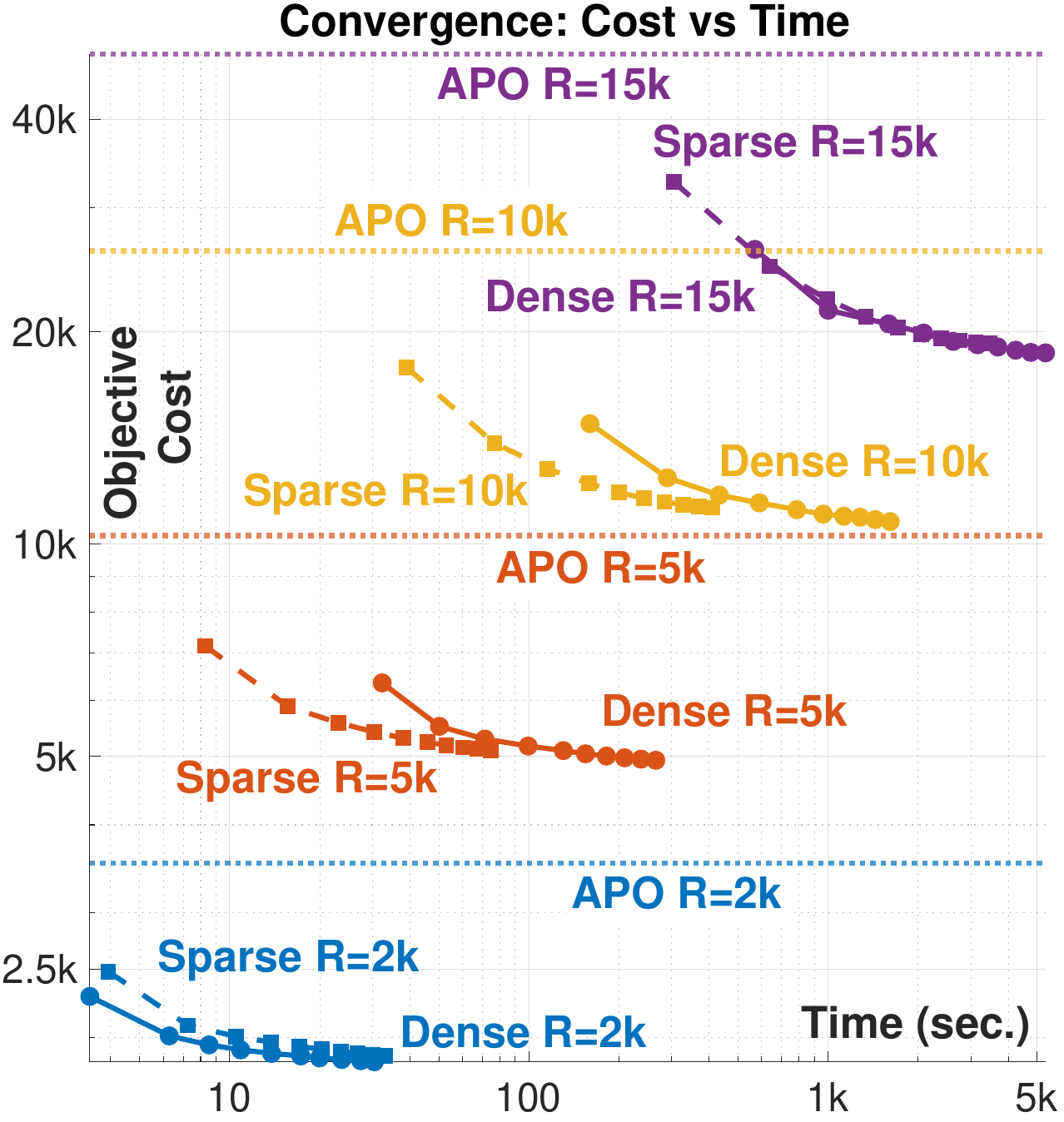}
\caption{\textbf{Left}: Performance comparison on Twitch ENGB dataset ($|\cE|\!=\!70648$).
\textbf{Right}: The convergence of cost vs time across $|\cR|$.
Here `k' means 1000.
}
\label{fig:exp5:all}
\end{figure}
\section{Conclusion}\label{sec:conc}
We identified Projected Performative Optimization ($\Pi$PO) as a new problem class that strictly generalizes standard PO by incorporating projection operators within the distribution map.
Applied to Optimal Network Pricing under the Oblivious User model, we showed that the projection induces Dead Zones, CQ failure, and PS collapse (phenomena absent from standard PO), and developed a theoretical and algorithmic framework to solve $\Pi$PO problems.
The Projected DDD map explicitly accounts for nonsmooth capacity gating, which causes standard gradient-based solvers to degenerate to trivial sub-optima.
We proved that the expected performative objective is strongly convex with a unique global $\Pi$PO, and developed a sparsity-aware TR-SQP solver with a principled Dead Zone escape mechanism justified by CQ failure theory.
The analysis reveals sensitivity to the Dead Zone, explaining why Jacobian treatment is necessary.
Experiments demonstrate significant welfare gains over existing methods, 
showing that our framework provides a principled foundation for AI systems learning robust pricing and resource-allocation decisions across networked systems, including the Internet and cloud computing, social media recommendation, smart grid load balancing, and autonomous traffic control.

More broadly, we extend performative prediction, a paradigm for learning in environments that react to a model's own decisions, to the nonsmooth, constrained optimization regimes that real-world systems impose.
We show that ignoring capacity gating yields a trivial solution; by contrast, our optimality concept, convergence theory, and scalable solver together enable good decisions to be learned under such feedback, advancing data-driven decision-making under performative, constrained dynamics.

\paragraph{Future works}
We list two meaningful future directions.
First, replacing the box $[\zeros,x_{\max}]$ in $\Pi_{[\zeros,x_{\max}]}$ with other convex sets in the DDD model.
Second, designing more ONP-specific efficient TR-SQP remains open: in \cref{sec:exp:subsec:Twitch}, we found that the QP solver is a computational bottleneck.
Replacing the QP solver with a structure-exploiting surrogate allows further scaling.
Furthermore, a sparsity-aware GPU HPC framework would enable more efficient computing for huge-scale networks.
More broadly, extending $\Pi$PO to other AI decision-making applications, such as AI inference infrastructure pricing (for ChatGPT, Claude, Gemini, etc.), online ad bidding under budget constraints, ride-hailing surge pricing, remains a promising direction.
\\~\\
\noindent\textbf{Acknowledgment.}
The authors acknowledge the IRIDIS High Performance Computing Facility and support at the University of Southampton. 
Li and Stein acknowledge the support from the EPSRC through the Turing AI Acceleration Fellowship on Citizen-Centric AI Systems (EP/V022067/1) and the Future Electric Vehicle Energy Networks supporting Renewables Grant
(EP/W005883/1).
\bibliographystyle{plain} 
\bibliography{ref}

\section{Appendix}
\etocsetnexttocdepth{subsection} 
\localtableofcontents 

\subsection{Proof of Proposition~\ref{prop:CA}}\label{pf:thm:CA}
\noindent
The total system $\text{cost}_{\text{tot}}(x)\!=\!\langle x, \text{cost}_{\cR}(x) \rangle$.
The route cost vector expands as $\text{cost}_{\cR}(x) \!=\! A(\Diag(c^{\text{coe}}) A^\top x \!+\!c^{\text{os}})$.
So $\text{cost}_{\text{tot}}(x) \!=\! \langle x, A \Diag(c^{\text{coe}}) A^\top x \rangle \!+\! \langle x, A c^{\text{os}} \rangle$.
Matching with $\tfrac{1}{2} \langle Qx, x \rangle \!-\! \langle s, x \rangle$ gives $Q \!=\! 2 A \Diag(c^{\text{coe}}) A^\top$ and $s \!=\! -Ac^{\text{os}}$.
Since $c^{\text{coe}}_e \!\geq\! 0$, $\Diag(c^{\text{coe}}) \!\succeq\! \zeros$ and  $Q \!\succeq\! \zeros$.

\subsection{Proof of Lemma~\ref{gradientHessian_derivation}}\label{sec:app:pf:lem:gradient_derivation}
\noindent
From \eqref{prob:ONP-nox-SAA},
$\nabla_p f_N(p)\!=\!\Avg[(\nabla_p x^i_p)^\top (Qx^i_p \!-\!s)]\!+\! \lambda p$.
The Clarke chain rule gives $\nabla_p x_p\!=\!JB$\footnote{With $y\!=\!Bp+\zeta$ and $x=\Pi_{[\zeros, x_{\max}]}(y)$: $\nabla_p x_p\!=\!\nabla_y x_p \cdot \nabla_p y\!=\!JB$.}, yielding the stated gradient.
Similarly, $\nabla_p c_N(p)\!=\!-\! K \Avg[J^i] B$.
Differentiating $\nabla_p f_N$ and dropping $\nabla_p J^i$\footnote{$\nabla_p J^i\!=\!0$ a.e.\ since the Jacobian \eqref{eq:Jac_subdiff} is constant except at the measure-zero projection boundaries.} gives the stated Hessian.
$\nabla^2_p c_{N}(p)\!\overset{\mathrm{a.s.}}{=}\!0$ follows similarly.

\subsection{Proof of Proposition~\ref{prop:concentration_combined}}\label{sec:app:subsec:pf:thm:concentration_combined}
We first introduce a useful lemma.
\begin{lemma}[ONP is Bounded and Lipschitz]\label{lem:onp_constants}
Let $x^i_p = \Pi_{[\zeros,x_{\max}]}(Bp\!+\!\zeta^i)$.
Almost surely, the per-sample derivatives of $f(p,\zeta^i)$ are
\[
\nabla_p f(p,\zeta^i)
\!=\! B^\top\!J^i(Qx^i_p\!-\!s)\!+\!\lambda p,
~
\nabla_p^2 f(p,\zeta^i)
\!=\!B^\top\! J^iQJ^i B \!+\! \lambda I,
~
\nabla_p^3 f(p,\zeta^i)
\!\overset{\mathrm{a.s.}}{=}\!\zeros,
\]
with the following uniform bounds and Lipschitz constants:
\begin{center}\footnotesize
\renewcommand{\arraystretch}{1}
\begin{tabular}{c|ll}
&\hspace{-6pt} \textbf{Bound ($M$)} &\hspace{-6pt} \textbf{Lipschitz ($L$)} \\
\hline
$f$
  &\hspace{-6pt} $M_f \!=\! \tfrac{1}{2}\|Q\|_2\|x_{\max}\|_2^2 \!+\! \|s\|_2\|x_{\max}\|_2$
  &\hspace{-6pt} $L_f \!=\! M_g \!+\! \lambda\|p_{\max}\|_2$ \\
$\nabla f$
  &\hspace{-6pt} $M_g \!=\! \|B\|_2(\|Q\|_2\|x_{\max}\|_2 \!+\! \|s\|_2)$
  &\hspace{-6pt} $L_g \!=\! M_H \!+\! \lambda$ \\
$\nabla^2 f$
  &\hspace{-6pt} $M_H \!=\! \|B\|_2^2\|Q\|_2$
  &\hspace{-6pt} $L_H \!\overset{a.s.}{=}\! 0$, $L_H^{\IE} \!<\! \infty$ (Lemma~\ref{lem:as_smoothing})
\end{tabular}
\end{center}
\noindent
Note that $(\ftrue, f_N)$ share $L_f$ globally, and  $(\nabla \ftrue, \nabla f_N)$ share $L_g$ cell-wise ($L_g$ is undefined on kinks).
For Hessian, $\HN \!\coloneqq\! \nabla^2 f_N$ is piecewise constant in $p$ (local $L_H \!=\! 0$), with $\cO(M_H)$ jumps at cell boundaries ($L_H \!=\! \infty$ on kinks), while $\HIE \!\coloneqq\! \nabla^2\ftrue$ is $C^\infty$ and $L_H^{\IE}$-Lipschitz globally.\footnote{
$M_f$ bounds $|f|$ by $\|Q\|_2\tfrac{\|x\|_2^2}{2} \!+\! \|s\|_2\|x\|_2 \!\leq\! M_f$ and $\|x\|_2 \leq \|x_{\max}\|_2$.
$M_g$ bounds $\|\nabla f\|$ via $\|J\|_2 \!\leq\! 1$ and $\|x\|_2 \!\leq\! \|x_{\max}\|_2$.
$M_H$ bounds $\|\nabla^2\! f\|_2$ via $\|J\|_2 \!\leq\! 1$.
Finiteness of $L_H^{\IE}$ follows from Gaussian smoothing (Lemma~\ref{lem:as_smoothing}).
}
\end{lemma}

We now prove Proposition~\ref{prop:concentration_combined}.
\paragraph{Step.0 (Setup)}
Let $\HN(p)\!\coloneqq\!\nabla^2 f_N(p)$ and $\HIE(p)\!\coloneqq\!\nabla^2\ftrue(p)$.
Fix $\epsilon\!>\!0$, $\delta\!\in\!(0,1)$.
Let $M,L,L_H^{\IE}$ as in Lemma~\ref{lem:onp_constants}.
Set the covering radii $\tau$ as 
$  \tau_f \coloneqq \tfrac{\epsilon}{4L_f},
  \tau_g \coloneqq \tfrac{\epsilon}{4L_g},
  \tau_H \coloneqq \tfrac{\epsilon}{2L_H^{\IE}}
$.
Set net sizes $n^{\bowtie}$ as
$ n^{\bowtie}_{\tau_q}$ where $ q\in\{f,g,H\}$.

\paragraph{Step.1 (Cover)}
Fix radius $\tau\!>\!0$, pick a finite $\tau$-net $\{p_1,\dots,p_{n^{\bowtie}_\tau}\}\!\subset\!\cP$ with
$\forall p\!\in\!\cP,\;\exists j:\!\|p\!-\!p_j\|\!\leq\!\tau$.
For $\nabla f_N$ and $\HN$, due to kinks, the net is additionally
\emph{cell-respecting}: that each $p\!\in\!\cP$ is matched to a $p_j$ in the \emph{same polyhedral cell}\footnote{Cell-respecting means that the net points avoid cell boundaries.
This is feasible since boundaries form a finite measure-zero union of hyperplanes.
}.

\paragraph{Step.2 (Bernstein at grid points)}
For a fixed single point $p_j$, any SAA average $q_N(p_j)\!=\!\Avg\![q(p_j,\zeta^i)]$ with $|q|\!\leq\!M_q$ satisfies
\[
  \IP\!\left(\bigl|q_N(p_j)-q_\IE(p_j)\bigr|\!>\!t\right)
  \leq
  2\exp\!\left(-N t^2/(2M_q^2+\tfrac{4}{3}M_q t)\right),
  ~~
  q\in\{f,g,H\}.
\]
Let $t\!=\!\tfrac{\epsilon}{2}$ and take a union bound over all $n^{\bowtie}_{\tau_q}$ net points gives
\[
\IP\text{(Failure$_q$)} 
=
n^{\bowtie}_{\tau_q} 
\cdot 
2 \exp\!\left(
    \tfrac{-N(\epsilon/2)^2}{2M_q^2 + \tfrac{4}{3}M_q(\epsilon/2)}
  \right).
\tag{Failure Probability}
\]

\paragraph{Step.3a (Gap-filling for $f$)}
Both $f_N, \ftrue$ are $L_f$-Lipschitz globally (Lemma~\ref{lem:onp_constants}).
For any $p$ with nearest net point $p_j$ ($\|p-p_j\|\!\leq\!\tau_f$):
\[
\begin{array}{rcl}
  |f_N(p)\!-\!\ftrue(p)|
  &\!\leq\!&
  \underbrace{\bigl|f_N(p)\!-\!f_N(p_j)\bigr|}_{\leq L_f\tau_f}
  +
  \underbrace{\bigl|f_N(p_j)\!-\!\ftrue(p_j)\bigr|}_{\leq \epsilon/2 \text{ (Step.2)}}
  +
  \underbrace{\bigl|\ftrue(p_j)\!-\!\ftrue(p)\bigr|}_{\leq L_f\tau_f}
  \\
  &\!=\!&
  \tfrac{\epsilon}{2}+2L_f\tau_f.
\end{array}
\]
With $\tau_f \!=\! \epsilon/(4L_f)$, this gives $\sup_{p\in\cP}\!|f_N(p)\!-\!\ftrue(p)|\!\leq\!\epsilon$.

\paragraph{Step.3b (Gap-filling for $\nabla f$)}
Each per-sample gradient is $L_g$-Lipschitz \emph{within each cell}
With the cell-respecting net, for each $p$, take $p_j$ in the \emph{same cell}, and then $\|\nabla f_N(p)\!-\!\nabla f_N(p_j)\|\!\leq\! L_g\tau_g$.
Since $\nabla\ftrue$ is globally $L_g$-Lipschitz, following the same triangular inequality used in Step.3a gives $\|\nabla f_N(p)\!-\!\nabla\ftrue(p)\| \!\leq\! \tfrac{\epsilon}{2}\!+\!2L_g\tau_g$.
With $\tau_g\!=\!\tfrac{\epsilon}{4L_g}$, this gives $\sup_{p\in\cP}\!\|\nabla f_N(p)\!-\!\nabla\ftrue(p)\|\!\leq\!\epsilon$.

\paragraph{Step.3c (Gap-filling for $\HN$)}
$\HN$ is piecewise constant in $p$ (local $L_H\!=\!0$), but jumps by $\cO(M_H)$ across cell boundaries\footnote{
Standard $\epsilon$-net argument fails here.
If $p$ and $p_j$ straddle a kink,  we have $\|\HN(p)-\HN(p_j)\|=\cO(M_H)$ error regardless of $\tau_H$:  the gap is not controlled by grid size at all.
}.
For each $p$, take $p_j$ in the \emph{same cell} with $\|p-p_j\|\!\leq\!\tau_H$, and  $\HN(p)\!=\!\HN(p_j)$ by $H_N$ being piecewise constant.
Combined with the $L_H^{\IE}$-Lipschitzness of $\HIE$ (Lemma~\ref{lem:onp_constants}):
\[
\begin{array}{l}
\hspace{-3pt} \|\HN(p)\!-\!\HIE(p)\|
  \leq \hspace{-3pt}
  \underbrace{\|\HN(p)\!-\!\HN(p_j)\|}_{=0 \text{ (piecewise constant)}}
  \hspace{-1pt}
  \!+\!
  \underbrace{\|\HN(p_j)\!-\!\HIE(p_j)\|}_{\leq \epsilon/2}
  \!+\!
  \underbrace{\|\HIE(p_j)\!-\!\HIE(p)\|}_{\leq L_H^{\IE}\tau_H}
\\
\hspace{73pt}
=  \frac{\epsilon}{2}+L_H^{\IE}\tau_H,
\end{array}
\]
so $\tau_H\!=\!\tfrac{\epsilon}{2L_H^{\IE}}$ 
gives $\sup_{p\in\cP}\|\HN(p)\!-\!\HIE(p)\|\!\leq\!\epsilon$.

\paragraph{Step.4 (Union bound and solving for $N$)}
Now we set each of the three events to fail with failure probability at most $\delta/3$.
Set $\IP\text{(Failure$_q$)} \leq \delta/3$ in Step.2 for each $q$ with bound $M_q$ and net size $n^{\bowtie}_{\tau_q}$, requiring
\[
  2n^{\bowtie}_{\tau_q}\exp\!\left(
    \tfrac{-N(\epsilon/2)^2}{2M_q^2 + \tfrac{4}{3}M_q(\epsilon/2)}
  \right) \leq \tfrac{\delta}{3}
~
\overset{\text{solve for }N}{\implies}~
  N_q^{\mathrm{unif}}.
\]
Now take $M\!=\!\max(M_f,M_g,M_H)$,
$N_{SAA}\!=\!\max(N_f^{\mathrm{unif}}, N_{g}^{\mathrm{unif}}, N_{H}^{\mathrm{unif}})$,
and a union bound over the three events (each at confidence $1-\delta/3$):
\[
  \IP\!\bigg(
    \sup_{p\in\cP}\bigl|f_N\!-\!\ftrue\bigr|\!>\!\epsilon
    ~\text{or}~
    \sup_{p\in\cP}\bigl\|\nabla f_N\!-\!\nabla\ftrue\bigr\|\!>\!\epsilon
    ~\text{or}~
    \sup_{p\in\cP}\bigl\|\HN\!-\!\HIE\bigr\|\!>\!\epsilon
  \bigg)
~\leq~
\delta.
\]
This gives\footnote{Bernstein is applied component-wise ($|\cR|$ for the gradient,
$|\cR|^2$ for the Hessian); the union bounds contribute
$\ln|\cR|$ and $2\ln|\cR|$ terms absorbed into the final logarithm.
}
\[
\begin{array}{l}
N_{SAA}
=
\max\Bigg\{\!
\underbrace{\cO\!\left(
\tfrac{M_f^2}{\epsilon^2}
\ln\tfrac{n^{\bowtie}_{\tau_f}}{\delta}
\right)}_{N_f^{\mathrm{unif}}}
,
\underbrace{
\cO\!\left(
\tfrac{M_g^2}{\epsilon^2}
\ln\tfrac{n^{\bowtie}_{\tau_g}|\cR|}{\delta}
\right)}_{N_g^{\mathrm{unif}}},
\underbrace{
\cO\!\left(
\tfrac{M_H^2}{\epsilon^2}
\ln\tfrac{n^{\bowtie}_{\tau_H}|\cR|^2}{\delta}
\right)}_{N_H^{\mathrm{unif}}}\!
\Bigg\}
\\[3pt]
\hspace{0mm}\leq
\cO\!\Big(
\tfrac{M^2}{\epsilon^2}
\ln\tfrac{
\max\{
n^{\bowtie}_{\tau_f},
n^{\bowtie}_{\tau_g},
n^{\bowtie}_{\tau_H}
\}
|\cR|^2
}{\delta}
\Big)
\overset{n^{\bowtie}_\epsilon
\coloneqq
\max\{
    n^{\bowtie}_{\tau_f},
    n^{\bowtie}_{\tau_g},
    n^{\bowtie}_{\tau_H}
    \}}{=}
\cO\!\left(
\tfrac{M^2}{\epsilon^2}
\ln\tfrac{
n^{\bowtie}_{\epsilon}
|\cR|
}{\delta}
\right),
\end{array}
\]
where the largest covering number corresponds to the smallest radius\footnote{
Covering numbers are monotone in the radius
($
\tau_1 \leq \tau_2
\Longrightarrow
n^{\bowtie}_{\tau_1}
\geq
n^{\bowtie}_{\tau_2}
$).}:
\[
n^{\bowtie}_\epsilon
=
n^{\bowtie}_{\min\{
\tau_f,\tau_g,\tau_H
\}}
=
n^{\bowtie}_{
\min\{
\epsilon/4L_f, ~ \epsilon/4L_g, ~\epsilon/2L_H^{\IE}
\}
}.
\tag{\#}
\]
\paragraph{Step.5 (Explicit $N_{SAA}$)}
For $\cP\!\subset\!\IR^{|\cR|}$ with diameter $D\!\coloneqq\!\|p_{\max}\!-\!p_{\min}\|_2$, 
 we have the covering number
$n^{\bowtie}_{\tau} \leq \left(\tfrac{2D}{\tau}+1\right)^{|\cR|}$ \cite{vershynin2018high}.
Following (\#), we get
\[
\ln n^{\bowtie}_{\epsilon}
\!\leq\!
\displaystyle
|\cR|\ln\!\left(\tfrac{2D}{
\min\left\{
\epsilon/4L_f, 
\epsilon/4L_g, 
\epsilon/2L_H^{\IE}
\right\}
}
\!+\!1\right)
\!=\!
\displaystyle
|\cR|\ln\!\left(\tfrac{D}{\epsilon
\min\{
1/2L_f, 
1/2L_g, 
1/L_H^{\IE}
\}
}
\!+\!1\right).
\]
Hence $\ln n^{\bowtie}_{\epsilon}
=
\cO\!\left(|\cR|\ln\tfrac{DL
}{\epsilon}\right)$, where $L=\max(2L_f, 2L_g,L_H^{\ftrue})$. 
Finally, integrating $\ln n_\epsilon^{\bowtie}$ into $N_{SAA}$ gives
\[
\boxed{
    N_{SAA}
    =
    \cO\!\Big(
      M^2\epsilon^{-2}
      \!\left[|\cR|\ln\tfrac{\|p_{\max}\!-\!p_{\min}\|_2 L}{\epsilon}+\ln\tfrac{|\cR|}{\delta}\right]
\Big).}
\]

\paragraph{Error Gap}
Let $p_N^*\coloneqq\arg\min_{p\in\cP}f_N(p)$ and
$p_\IE^*\coloneqq\arg\min_{p\in\cP}\ftrue(p)$.
From the uniform bound above, $\sup_{p\in\cP}|f_N(p)-\ftrue(p)|\leq\epsilon$
w.p.  $\ge1-\delta$.
\\
Upper bound: $\min_{p} f_N(p)\!=\!f_N(p_N^*) \geq \ftrue(p_N^*) \!-\!\epsilon \geq \min_{p} \ftrue(p) \!-\!\epsilon$.
\\
Lower bound: $\min_{p} f_N(p) \leq f_N(p_\IE^*) \leq \ftrue(p_\IE^*) \!+\!\epsilon\!=\!\min_{p} \ftrue(p) \!+\!\epsilon$.
\\
Combining gives $\left|\min_{p\in\cP}f_N(p)\;-\;\min_{p\in\cP}\ftrue(p)\right|\;\leq\epsilon$.

\subsection{Proof of Prop.~\ref{prop:ftrue_no_DZ_min}}
\noindent
Let $p^\bullet_N \neq \zeros$ be a stationary point of $f_N$, we first show by contradiction that $p^\bullet_N \notin \mathrm{int}\DZ(f_N)$.
Assume $\zeros \neq p^\bullet_N \in \mathrm{int}\DZ(f_N)$ is an interior stationary point.
\begin{enumerate}[leftmargin=*, label=(\roman*)]\setlength{\itemsep}{-2pt}
\item By \eqref{eq:DZ_grad_collapse}, $\nabla_p f_N(p^\bullet_N) = \lambda p^\bullet_N$.
\item Stationarity gives $\lambda p^\bullet_N = \zeros$. 
\item (i), (ii) give $\lambda p^\bullet_N\!=\!\zeros \overset{\lambda > 0}{\implies}$ $p^\bullet_N\!=\!\zeros$, contradiction.
\end{enumerate}
As we proved no nonzero stationary point in DZ, thus all nonzero minimizers cannot lie in DZ. 

It remains to exclude the zero minimizer of $f_N$ in DZ.
By Prop.~\ref{prop:PS_eq_PO},  for $\ftrue$, PS ($p^*_{\IE}\!=\!\zeros$) $\,=\!\Pi$PO holds only in degenerate cases. 
For the surrogate $f_N$ of $\ftrue$ that excludes PS, the minimizer $p^*_N$ 
lies outside $\DZ(f_N)$.

\subsection{$\pi_r$ of $\ftrue$ and $\ftrue$ has no dead zone}\label{sec:app:pf:prop:DZfE_no_PO_pt_pi_r_positive}
\noindent
By the Jacobian cases and Lemma~\ref{lem:as_smoothing}:
\[
\pi_r\!=\!
\IE[J_{rr}]
= \int_{0 < (Bp+\zeta)_r < (x_{\max})_r} \hspace{-5mm}\phi(\zeta) d\zeta
= \IP(0 < (Bp+\zeta)_r < (x_{\max})_r).
\]
For Gaussian $\zeta$:
$\IE[J_{rr}] =
\Phi\bigl( \tfrac{(x_{\max})_r -(Bp)_r -\mu_r}{\sigma_r} \bigr)
\!-\!\Phi\bigl( -\tfrac{(Bp)_r +\mu_r}{\sigma_r} \bigr)$,
which is strictly positive when $(x_{\max})_r > 0$ since $\Phi$ is strictly increasing.
So we proved $\pi_r(p)>0$ for \textbf{all} $r\in\cR$ and \textbf{all} feasible $p$.
This proved that $\DZ(\ftrue)\!=\!\varnothing$.

\subsection{Proof of Corollary~\ref{cor:pIE_not_in_DZ}}
We prove the Corollary in 6 steps.
\begin{enumerate}[leftmargin=*]\setlength{\itemsep}{0pt}
\item Let $\epsilon_1 \coloneqq \dist(p^*_N, \DZ(f_N))$.
By Prop.~\ref{prop:ftrue_no_DZ_min} we have $\epsilon_1 >0$.

\item Prop.~\ref{prop:concentration_combined} gives $|f_N(p^*_N) \!-\!\ftrue(p^*_{\IE})| \leq \epsilon$ w.p.\ $1-\delta$.

\item Now will show that $(\ast) : \tfrac{\alpha}{2}\|p^*_{\IE} \!-\!p^*_N\|_2^2 
\leq \ftrue(p^*_N) \!-\!\ftrue(p^*_{\IE}) \leq 2\epsilon$
\begin{itemize}[leftmargin=*]\setlength{\itemsep}{0pt}
    \item We have  
    $\ftrue(p^*_N) 
\geq  \ftrue(p^*_{\IE}) 
+
\big\langle
\nabla \ftrue(p^*_{\IE}), ~ p^*_N \!-\!p^*_{\IE}
\big\rangle
+
\tfrac{\alpha}{2}\|p^*_N \!-\!p^*_{\IE}\|_2^2$ since $\ftrue$ is $\alpha$-strongly convex (Theorem~\ref{thm:fE_strong_cvx}).
Put $\nabla\ftrue(p^*_{\IE})\!=\!\zeros$ gives
\[
\begin{array}{rcll}
\tfrac{\alpha}{2}\|p^*_N \!-\!p^*_{\IE}\|_2^2
&\leq&
\ftrue(p^*_N)  \!-\!\ftrue(p^*_{\IE}) 
\!=\!
\ftrue(p^*_N)\!-\!f_N(p^*_N) \!+\!f_N(p^*_N)  \!-\!\ftrue(p^*_{\IE}) 
\\
&\leq&
\big|  \ftrue(p^*_N) \!-\!f_N(p^*_N) \big| 
\!+\!
\big| f_N(p^*_N)  \!-\!\ftrue(p^*_{\IE})  \big|
~~ \text{(tri. iq.)}
\\[2pt]
&\leq&\displaystyle
\big|  \min_{p} f_N \!-\!\ftrue(p^*_N) \big| 
\!+\!
\sup_p
\big| f_N  \!-\!\ftrue  \big|
\\
&\leq&
\hspace{-0.33cm}
\underbrace{| \min_{p}f_N \!-\!\min_{p}\ftrue  |}_{\leq\,\epsilon \text{ by Error gap in Prop.}\ref{prop:concentration_combined}} 
\!+\!\underbrace{\sup_{p}|f_N \!-\!\ftrue|}_{\leq\,\epsilon \text{ by Prop.}\ref{prop:concentration_combined}} 
=~ 2\epsilon ~ \text{w.p. }1-\delta
\end{array}
\]
\end{itemize}
\item ($\ast$) gives $\|p^*_{\IE} \!-\!p^*_N\|_2 \leq \sqrt{4\epsilon/\alpha}$.
Pick $\epsilon < \alpha\epsilon_1^2/4$  ensures $\|p^*_{\IE} \!-\!p^*_N\|_2 < \epsilon_1$.

\item Combine the results in step-1 and step-4 give
$
\|p^*_{\IE} \!-\!p^*_N\|_2 
\leq \epsilon_1
$.

\item By the reverse triangle inequality\footnote{Take $s \in \DZ(f_N)$. 
Triangle inequality gives $|p^*_N \!-\!s\| \;\leq\; \|p^*_N \!-\!p^*_{\IE}\| \!+\!\|p^*_{\IE} \!-\!s\|$.
Take the infimum over all $s \in \DZ(f_N)$ on both sides:
\[
\underbrace{\inf_{s \in \DZ(f_N)} \|p^*_N \!-\!s\|}_{\eqqcolon \dist(p^*_N, \DZ(f_N))}
\;\leq\;
\|p^*_N \!-\!p^*_{\IE}\|
\;+\;
\underbrace{
\inf_{s \in \DZ(f_N)} \|p^*_{\IE} \!-\!s\|}_{\eqqcolon \dist(p^*_{\IE}, \DZ(f_N))}.
\]
Hence $\dist(p^*_N, \DZ(f_N))
\;\leq\;
\|p^*_N \!-\!p^*_{\IE}\|
\;+\;
\dist(p^*_{\IE}, \DZ(f_N)).$
Subtract $\|p^*_N \!-\!p^*_{\IE}\|$ from both sides, and swapping the sides gives exactly the inequality.},
\[
\dist(p^*_{\IE}, \DZ(f_N))
\geq \dist(p^*_N, \DZ(f_N)) \!-\!\|p^*_{\IE} \!-\!p^*_N\|_2
= \epsilon_1 \!-\!\|p^*_{\IE} \!-\!p^*_N\|_2 > 0,
\]
where the last inequality uses step-4.
Thus $\dist(p^*_{\IE}, \DZ(f_N)) > 0$ w.p. $1-\delta$, which implies
$\IP(p^*_{\IE} \in \DZ(f_N)) \leq \delta$.
\end{enumerate}

\subsection{Proof of Proposition~\ref{prop:fail_CQ}}\noindent
As $p_r \in \DZ_r(f_N)$, $J^i_{rr}=0$ for all $i$, so $[\nabla_p c_N]_{r,\cdot}=\zeros$ by Lemma~\ref{gradientHessian_derivation}. LICQ fails (zero vectors are linearly dependent) and MFCQ fails ($\zeros d=0\not<0$ for all $d$).

\subsection{Proof of Theorem~\ref{thm:fE_strong_cvx}}\noindent
(For compactness we ignore the additive constant $\lambda$.)

\noindent
\textbf{Part (i): PD Hessian and Strong convexity.}
Let $v \in \IR^{|\cR|}$ and $u\!=\!Bv$.
Let $\omega\!=\!\Diag(J) \in \{0,1\}^{|\cR|}$ and $\pi\!=\!\IE[\omega]$.
Then $Ju\!=\!\omega \odot u$ and
$v^\top \IE[H_N] v\!=\!u^\top (\IE[\omega \omega^\top] \odot Q) u$.
Writing $\IE[\omega \omega^\top]\!=\!\Sigma_\omega \!+\!\pi\pi^\top$ with $\Sigma_\omega \succeq \zeros$, and applying the Schur Product Theorem ($\Sigma_\omega \odot Q \succeq \zeros$):
$\IE[H_N] \succeq B^\top (\pi\pi^\top \odot Q) B\!=\!B^\top \Diag(\pi) Q \Diag(\pi) B$.
Since $\Diag(\pi) \succeq \pi_{\min} I$:
\[
u^\top \Diag(\pi) Q \Diag(\pi) u
= \| Q^{1/2} \Diag(\pi) u \|^2
\geq \pi_{\min}^2 \lambda_{\min}(Q) \sigma_{\min}^2(B) \|v\|^2,
\]
giving $\IE[H_N] \succeq \pi_{\min}^2\lambda_{\min}(Q)\sigma_{\min}^2(B) I$.
Now, $\nabla^2 \ftrue(p) \succeq \alpha I \succ \zeros$ means $\ftrue$ is $\alpha$-strong convex.

\noindent \textbf{Part (ii): CQ holds.}
By \cref{sec:app:pf:prop:DZfE_no_PO_pt_pi_r_positive},
$\pi_r(p) \!>\! 0$ for all $r \!\in\! \cR$ and all feasible $p$,
so $\Diag(\pi) \!\succ\! \zeros$ and $\nabla_p c_{\IE}(p) \!=\! -KB\Diag(\pi)$ is full rank from (A3).

Hence, MFCQ holds for $\ftrue$,
making KKT necessary; sufficiency follows from $\alpha$-strong convexity.
Thus, \eqref{eq:KKT_true} characterizes the unique global 
$p^*_{\Pi\text{PO}}$.

\subsection{Proof of Proposition~\ref{prop:reg_amp_gap}}\label{sec:app:pf:prop:reg_amp_gap}\noindent
We analyze the case that $p^*_{\Pi\text{PO}}$, $p^*_{APO}$ are interior
(taking $\gamma\!=\!\zeros$ for exposition).\footnote{Boundary solutions are
viewed as limits of $\epsilon$-relaxed interior problems.}

\paragraph{Part (i): Upper regime}
By definition, $\hat{x}_{APO}\!=\!y \geq \Pi_{x_{\max}}(y)\!=\!x_{\Pi\text{PO}}$
(APO overflows).
$-B^\top Q \overset{A4}{\geq} \zeros$ element-wise, so multiplying by
$-\tfrac{1}{\lambda}B^\top Q$ preserves direction: $p^*_{APO}
\!=\!-\tfrac{1}{\lambda}B^\top \IE[Q\hat{x}^{*}_{APO}-s]
\geq -\tfrac{1}{\lambda}B^\top \IE[J(Qx^{*}_{\Pi\text{PO}}-s)]
\!=\!p^*_{\Pi\text{PO}}.$

For $\texttt{Gap}$: since $x \in [\zeros, x_{\max}]$, let
$g(p)\!=\!\IE[\langle\tfrac{1}{2}Qx_p \!-\!s, x_p\rangle]$ be bounded by $M$.
Then $\texttt{Gap} \geq -2M \!+\!\tfrac{\lambda}{2}(\|p_{APO}\|^2 -
\|p_{\Pi\text{PO}}\|^2)$.
Since $p_{APO} \geq p_{\Pi\text{PO}}$, the price term is positive and dominates
for large $\lambda$.

\paragraph{Part (ii): Lower regime}
We use two helper lemmas.

\begin{lemma}\label{lem:monotone_obj_p}
For sufficiently small $\lambda$, a smaller $p$ increases $\ftrue$.
\end{lemma}
\begin{proof}
$\ftrue(p) \approx \IE[\langle \tfrac{1}{2}Qx \!-\!s, x \rangle]$ for
small $\lambda$.
Since $\nabla_p \IE[x]\!=\!\IE[JB] \preceq \zeros$ a.s., lower prices
increase expected flow, raising $\ftrue$.
\end{proof}

\begin{lemma}[Restricted Feasible Sets]\label{lem:optgap_PO_diffx}
For $\cX_{\min}\!=\!\{x \in \cX : x \ge x_{\min}\} \subset \cX$:
$f(p_{\Pi\text{PO}}^*(\cX_{\min})) \geq f(p_{\Pi\text{PO}}^*(\cX))$.
\end{lemma}

\noindent
By definition, $\hat{x}_{APO}\!=\!y \leq \Pi_{x_{\min}}(y)\!=\!x_{\Pi\text{PO}}$
(APO underflows).
By (A4), mirroring Part~(i):
$p_{APO}^* \leq -\tfrac{1}{\lambda}B^\top \IE[Qx_{\Pi\text{PO}}^{*}-s]$.
Since $\cX_{\min} \subset \cX$ forces higher prices,
$p_{APO}^*(\cX) \leq p_{APO}^*(\cX_{\min})$, and:
\[
\ftrue(p_{APO}^*(\cX)) \geq \ftrue(p_{APO}^*(\cX_{\min}))
= \ftrue(p_{\Pi\text{PO}}^*(\cX_{\min}))
\geq \ftrue(p_{\Pi\text{PO}}^*(\cX)),
\]
by Lemma~\ref{lem:monotone_obj_p} and
Lemma~\ref{lem:optgap_PO_diffx}.

\subsection{Proof of Lemma~\ref{lem:grad_mismatch}}\noindent
APO stationarity: $\IE[B^\top(Qx \!-\!s)] \!+\!\lambda p_{APO}\!=\!\zeros$.
Substituting into the true PO gradient:
$\nabla \ftrue(p_{APO})\!=\!\IE[B^\top J(Qx \!-\!s)] \!+\!\lambda p_{APO}\!=\!B^\top \IE[(J \!-\!I)(Qx \!-\!s)]$.

Since $\ftrue$ is quadratic (Theorem~\ref{thm:fE_strong_cvx}), Taylor expansion around $p_{\Pi\text{PO}}$ (where $\nabla\ftrue\!=\!\zeros$) gives
$\ftrue(p_{APO}) \!-\!\ftrue(p_{\Pi\text{PO}})\!=\!\tfrac{1}{2} \Delta p^\top H \Delta p$
with $H\!=\!\nabla_p^2 \ftrue(p_{\Pi\text{PO}}) \succ \zeros$ and $\Delta p\!=\!H^{-1}\nabla \ftrue(p_{APO})$.

Substituting: $\ftrue(p_{APO}) \!-\!\ftrue(p_{\Pi\text{PO}})\!=\!\tfrac{1}{2} \nabla\ftrue(p_{APO})^\top H^{-1} \nabla\ftrue(p_{APO})$.
The Rayleigh quotient gives $\alpha^{-1} \in [\lambda_{\max}(H)^{-1}, \lambda_{\min}(H)^{-1}]$, establishing the result.

\subsection{Proof of Theorem~\ref{thm:PS_failure}}\label{app:pf:thm:PS_failure}\noindent
By PS definition,
$p_{PS}\!=\!\argmin_{p \in [p_{\min},p_{\max}]}
  \IE_{x \sim \cD(p_{PS})}[\text{cost}_{\text{tot}}(x)] \!+\!\tfrac{\lambda}{2}\|p\|_2^2$.
At the PS equilibrium, the expectation is over the distribution \emph{fixed} at $\cD(p_{PS})$, so
$C \coloneqq \IE_{x \sim \cD(p_{PS})}[\tfrac{1}{2}\langle Qx-s,x\rangle]$
is constant in $p$ and $p_{PS}\!=\!\displaystyle \argmin_{p \in [p_{\min},p_{\max}]} \tfrac{\lambda}{2}\|p\|_2^2$.

\subsection{Proof of Lemma~\ref{lem:finite_dz_escape}}\noindent
\textbf{Case (i), $\|\mathrm{res}_t\|_\infty > 0$:}
Let $\Delta_{\mathrm{entry}}$ be the TR radius upon entering \DZ\ (frozen during DZ iterations).
The reflection step uses $d_t\!=\!-\alpha_t B^{-1}\mathrm{res}_t$ with $\alpha_t\!=\!\min(1, \Delta_{\mathrm{entry}} / \|B^{-1}\mathrm{res}_t\|_\infty)$.

The new residual satisfies $\|\mathrm{res}_{t+1}\|_1 \leq (1-\alpha_t)\|\mathrm{res}_t\|_1$.
If $\alpha_t\!=\!1$: escape in one step.
If $\alpha_t < 1$: each step decreases $\|\mathrm{res}\|_1$ by at least $c_{\mathrm{dec}}\!=\!\Delta_{\mathrm{entry}}/\|B^{-1}\|_\infty > 0$, so escape occurs in at most $T_{\DZ} \leq \lceil \|\mathrm{res}_0\|_1 / c_{\mathrm{dec}} \rceil$ steps.

\textit{Descent of $\Phi$:} since
\begin{itemize}\setlength{\itemsep}{-1pt}
    \item $\langle \nabla_p V, d_t \rangle\!=\!-\alpha_t\|\mathrm{res}_t\|_1 < 0$
    \item $\langle \nabla_p f_N, d_t \rangle \leq \alpha_t\lambda\|p_{\max}\|_2\|B^{-1}\|_\infty\|\mathrm{res}_t\|_1$,
\end{itemize}
we have $\Phi(p_{t+1}) \!-\!\Phi(p_t) \leq -\alpha_t(\mu_V \!-\!\lambda\|p_{\max}\|_2\|B^{-1}\|_\infty)\|\mathrm{res}_t\|_1\!=\!-c_\Phi \!<\! 0$.

\noindent
\textbf{Case (ii), $\|\mathrm{res}_t\|_\infty\!=\!0$:}
Under $\Sigma \!\succ\! \zeros$, this event has probability zero (finite-precision artifact).
The nudge step displaces $\bar{y}_r$ by $B_{rr}(\pm\varepsilon_0) \neq 0$, either exiting \DZ\ or producing $\|\mathrm{res}_{t+1}\|_\infty \!>\! 0$ (entering case (i)), in exactly one step.
$\Phi$ decreases by at least $\mu_V\varepsilon_0\min_r|B_{rr}| \!>\! 0$.

\subsection{Proof of Theorem~\ref{thm:global_convergence}}
\begin{lemma}[Penalty Stabilization and Sufficient Descent]\label{lem:penalty_bounded_and_Sufficient_Descent}
Assume $\pi_r \!>\! 0$ for at least one route per commodity (MFCQ holds).
Then $\nu_t$ stabilizes at a finite $\nu^*$ after finitely many iterations.
At $\nu^*$, for non-stationary $p_t$:
$\mathrm{Pred}_t(d_t) \geq \kappa \min(\dist(\zeros, \partial \phi(p_t, \nu^*)), \Delta_t)$.
\end{lemma}
\begin{proof}
Total DZ steps are at most $(\Phi_0 \!-\! \Phi_{\min})/c_\Phi \!<\! \infty$ (Lemma~\ref{lem:finite_dz_escape}).
For large $t$, iterates stay in the active region where MFCQ holds and dual multipliers are uniformly bounded, so $\nu_t$ stabilizes.
Sufficient descent follows from standard SQP theory \cite{nocedal2006numerical} in the active region, where Clarke subdifferential reduces to the classical gradient.
\end{proof}

\noindent\textit{Proof of theorem.}
For large $t, \{p_t\}$ consists of active-region steps
(finite total DZ steps by Lemma~\ref{lem:finite_dz_escape} and 
$\mu_V\!>\!\lambda\|p_{\max}\|_2\|B^{-1}\|_\infty$ by assumption).
$\phi$ is bounded below:\!
$\sum_{t=K}^{\infty} \mathrm{Ared}_t \!\leq\! \phi_K \!-\! \phi_{\min} \!<\! \infty$,
so $\mathrm{Ared}_t\!\to\!0$.
On successful steps, $\mathrm{Ared}_t \!\geq\! \eta_1\mathrm{Pred}_t$, so $\mathrm{Pred}_t \!\to\! 0$.
By Lemma~\ref{lem:penalty_bounded_and_Sufficient_Descent},  we have 
$\mathrm{Pred}_t \!\geq\! \kappa\min(\dist(\zeros,\partial\phi(p_t,\nu^*)),\Delta_t)$.
If $\dist(\zeros,\partial\phi(p_t,\nu^*))\!\not\to\!0$, then 
$\mathrm{Pred}_t \!\geq\! \kappa\min(c,\Delta_t)$ for some $c\!>\!0$.
But $\mathrm{Pred}_t\to 0$ forces $\Delta_t\to 0$, 
contradicting the TR expansion rule on repeatedly successful steps.
Hence $\dist(\zeros,\partial\phi(p_t,\nu^*))\to 0$. \qed

\subsection{Proof of Corollary~\ref{cor:end_to_end}}
\noindent
\textbf{1.\,Finite DZ escape.}
Lemma~\ref{lem:finite_dz_escape} relocates to $p'$ with $\pi_{\min}(p')\!>\!0$ in at most $(\Phi_0-\Phi_{\min})/c_\Phi$ iterations.

\noindent
\textbf{2.\,Clarke-stationarity.}
$\|\nabla f_N(\hat{p}_N^*)\|\!\overset{a.s.}{\leq}\! \epsilon_{opt}$ (Theorem\,\ref{thm:global_convergence}) in the active region.

\noindent
\textbf{3.\,Gradient transfer.}
By Prop.~\ref{prop:concentration_combined}: $\|\nabla f_N(\hat{p}_N^*) \!-\!\nabla \ftrue(\hat{p}_N^*)\| \leq \epsilon$ w.p.\ $1-\delta$.
Triangle inequality gives $\|\nabla \ftrue(\hat{p}_N^*)\| \leq \epsilon \!+\!\epsilon_{opt}$.

\noindent 
\textbf{4.\,Distance bound.}
By $\alpha$-strong convexity(Theorem~\ref{thm:fE_strong_cvx}):
$\alpha\|\hat{p}_N^* \!-\!p^*_{\IE}\| \leq \|\nabla \ftrue(\hat{p}_N^*)\| \leq \epsilon \!+\!\epsilon_{opt}$.

\subsection{Details of Exp1}\label{app:exp_setup}
\paragraph{Graph}
$\cV\!=\!\{s,A,B,t\}$, $\cE\!=\!\{sA, At, AB, sB, Bt\}$, $\cR\!=\!\{r_1{:}sAt, r_2{:}sABt, r_3{:}sBt\}$ (3 routes), $K\!=\!\ones_{1 \times 3}$.
Edge costs follow \cite{braess1968paradoxon}: $sA, Bt$ have congestion-sensitive latency (equal to flow); $sB, At$ have constant latency $1$; bridge $AB$ has zero cost.

\begin{figure}[h!]
\label{fig:Braess_graph_details}
\centering\footnotesize
\setlength{\tabcolsep}{2pt}
\renewcommand{\arraystretch}{0.85}
\begin{tabular}{@{}c c@{}}
\begin{tabular}{l|ccccc}
route & $sA$ & $At$ & $AB$ & $sB$ & $Bt$ \\
\hline
$r_1$       & 1 & 1 & 0 & 0 & 0 \\
$r_2$ & 1 & 0 & 1 & 0 & 1 \\
$r_3$       & 0 & 0 & 0 & 1 & 1
\end{tabular}
~&
\begin{tabular}{l|cc}
edge & $c^{\text{coe}}_e$ & $c^{\text{os}}_e$ \\
\hline
$sA$ & 1 & 0 \\
$At$ & 0 & 1 \\
$AB$ & 0 & 0 \\
$sB$ & 0 & 1 \\
$Bt$ & 1 & 0
\end{tabular}
\end{tabular}
\caption{Braess network: assignment matrix $A$ and edge costs.}
\end{figure}

\paragraph{Baselines}
\textbf{RR:} PS step solved exactly by Interior Point to ensure failure is theoretical, not numerical.
\textbf{GD/SGD:} Penalized objective with projected box constraints; SGD uses mini-batches from fixed $\{\zeta_i\}$.
\textbf{APO:} Ignores $\Pi$; uses $\nabla \tilde{f}\!=\!B^\top(Qx \!-\!s) \!+\!\lambda p$ and $H_{\text{GGN}}$.

\paragraph{Parameters}
$\lambda\!=\!10^{-3}$, $N=3000$, $x_{\max}\!=\!l\!=\!\ones$, $B\!=\!-2I$, $x_{\min}\!=\!p_{\min}\!=\!\zeros$, $p_{\max}\!=\!10 \cdot \ones$.
$\mu\!=\!10 \cdot \ones$, $\sigma\!=\!0.05$, $\rho\!=\!0.01$ (deeply congested).
$p_0\!=\!B^{\dagger}(\tfrac{x_{\max}}{2}-\mu)$.

\subsection{Details of Exp2}\label{app:exp_setup_exp2}
\paragraph{Parameters}
For the second row of Fig~\ref{fig:exp2traj}, $B\!=\!\tfrac{1}{2}
\Bigl[ \begin{smallmatrix}
-7 & 1 & 5 \\
 3 &-8 & 4 \\
 2 & 6 &-9
\end{smallmatrix} \Bigr]
$.

\subsection{Details of Exp3}\label{app:exp_setup_exp3}
\paragraph{Parameters}
Braess network with $\mu\!=\!5 \cdot \ones$, $\sigma\!=\!2$, $\rho\!=\!0.01$, $l\!=\!\ones$, $N=3000$, $p_{\min}\!=\!\zeros$, $p_{\max}\!=\!10\cdot\ones$, $B\!=\!-2I$.

\paragraph{Evaluation}
Sweep $\lambda \in [10^{-2}, 10^2]$ and $x_{\max} \in [0.3, 7]$, averaging over $25$ trials.
Final expected cost evaluated with $10^4$ Monte Carlo samples from $\cN(\mu, \Sigma)$.
Infeasible APO trials (constraint violated or diverged) are shaded gray in Fig.\,\ref{fig:Exp3}.

\subsection{Sparsity-aware Efficient Computation}\label{sec:app:subsec:sparse}
\paragraph{$f_N$ on edge}
Rewriting $\langle Qx, x \rangle$ via edge flows $f_e\!=\!A^\top x$, we have that 
$f_N(p)\!=\!\Avg[ \sum_{e \in \cE} ( c^{\text{coe}}_e (f^{(i)}_e)^2 \!+\!c^{\text{os}}_e f^{(i)}_e )] \!+\!\tfrac{\lambda}{2}\|p\|_2^2$.
Cost $\cO(N|\cR|^2)$ drops to $\cO(N \nnz(A))$.

\paragraph{Gradient $\nabla_p f_N$}
Let $X_p\!=\![x^1_p, \dots, x^N_p] \in \IR^{|\cR| \times N}$ and $\Omega \!\in\! \{0,1\}^{|\cR| \times N}$ the saturation mask. 
Then $\nabla_p f_N(p)\!=\!\tfrac{1}{N} B^\top ( \Omega \odot (Q X_p \!-\! s \ones_N^\top) ) \ones_N \!+\!\lambda p$. The elasticity matrix $B$ is block-sparse: each commodity admits at most $k$ alternative routes, giving support $\cB\!=\!\{(r,r') : B_{r,r'} \neq 0\}$ with $|\cB|\!=\!|\cR| \!+\! |\cK|k(k{-}1) \ll |\cR|^2$ \footnote{Each commodity forms a $k \times k$ dense block with $k(k{-}1)$ off-diagonal entries. As $|\cK| k\!=\!|\cR|$, we have $|\cB|\!=\!\cO(k|\cR|)$.}. 
The parallel active-set loop has cost $\cO(\sum_{i=1}^N ( 2 \nnz(A_{S_i}) \!+\! |S_i| ) \!+\! |\cB|)$: Given $A$ (scaled by $\sqrt{c^{\text{coe}}}$), $B, s$, compute $X_p, \Omega$, set $M\!=\!\zeros^{|\cR| \times N}$.
\\
\textbf{Parallel loop:} For each $i\!=\!1, \dots, N$:
\begin{enumerate}[leftmargin=*]\setlength{\itemsep}{-2pt}
    \item $S_i\!=\!\{ r : \Omega_{r,i}\!=\!1 \}$ \hfill active set
    \item $A_{S_i}\!=\!A(S_i, :)$, $J_i\!=\!\{ e : \exists r \in S_i,\, A_{r,e} \neq 0 \}$ \hfill active rows and edges
    \item $u_{J_i}\!=\!(A_{S_i}(:, J_i))^\top X_p(S_i, i)$ \hfill intermediate vector
    \item $M(S_i, i)\!=\!A_{S_i}(:, J_i)\, u_{J_i} \!-\!s(S_i)$ \hfill weighted residual
\end{enumerate}
Output: $m\!=\!\tfrac{1}{N} M \ones_N$;\quad $\nabla_p f_N(p)\!=\!\lambda p \!+\!B^\top m$ \hfill $\cO(|\cB|)$ sparse multiply

\paragraph{Gradient $\nabla_p c_N$}
Compute $d\!=\!\Omega \ones_N$ in $\cO(N |\cR|)$, and $\nabla_p c_N(p)\!=\!-\! \tfrac{1}{N} K \Diag(d)\, B$ in $\cO(\nnz(K) \cdot k)$, exploiting the block-sparse structure of $B$ \footnote{$K\Diag(d)$ costs $\cO(\nnz(K))$; the product with sparse $B$ costs $\cO(\nnz(K) \cdot k)$.}.

\paragraph{Hessian $\nabla_p^2 f_N$}
As $J^i Q J^i\!=\!(\omega^i {\omega^i}^\top)\odot Q$\footnote{$(J^iQJ^i)_{kl}=\omega^i_k Q_{kl} \omega^i_l$.} and $\sum_i \omega^i {\omega^i}^\top\!=\!\Omega\Omega^\top$:
$\nabla^2_p f_N (p)\!=\!\lambda I \!+\! \tfrac{1}{N} B^\top ( (\Omega\Omega^\top) \odot Q )  B$.
The active-set loop accumulates $Q_{S_i,S_i}$. The cost is $\cO(\sum_i |S_i|^2 (|J_i| \!+\!k))$\footnote{The loop costs $\cO(\sum_i |S_i|^2 |J_i|)$; $B^\top H_{\text{in}} B$ costs $\cO(|\cB| \cdot \bar{h}) \leq \cO(k \sum_i |S_i|^2)$, where $\bar{h}$ is the average number of nonzeros per row of $H_{\text{in}}$, since $|\cB|\!=\!\cO(k|\cR|)$ and $\bar{h} \leq \tfrac{1}{|\cR|}\sum_i |S_i|^2$.}.
\textbf{Parallel loop:} For each $i$: (1) Get $S_i, A_{S_i}, J_i$ as above and (2) $H_{\text{in}}(S_i, S_i) \mathrel{+}\!=\!2 A_{S_i}(:, J_i) (A_{S_i}(:, J_i))^\top$.
After the loop output $\nabla^2_p f_N(p)\!=\!\lambda I \!+\!\tfrac{1}{N} B^\top H_{\text{in}} B$.

\begin{table}[h!]
\centering
\caption{Per-iteration complexity. $S_i (J_i$): active routes (edges); $k$: routes per commodity.}
\footnotesize
\begin{tabular}{l|ccc}
\hline
 & $f_N(p)$ & $\nabla_p f_N(p)$ & $\nabla_p^2 f_N(p)$ \\
\hline
Sparse & $\cO(N\nnz(A))$ & $\cO(N\nnz(A)\!+\!k|\cR|)$ & $\cO(\sum_i |S_i|^2 (|J_i| \!+\!k))$
\\
Dense & $\cO(N|\cR|^2)$ & $\cO(N|\cR|^2)$ & $\cO(N|\cR|^2 \!+\!|\cR|^3)$
\end{tabular}
\end{table}

\paragraph{Accelerating SQP via low-rank structure}
Using the factorization $Q\!=\!FMF^\top$ with $\rk(A) \leq |\cE|$:
$\nabla^2_p f_N \overset{\text{a.s.}}{=} \U \M \U^\top \!+\! \lambda I$,
where $U_i\!=\!B^\top J^i F \in \IR^{|\cR| \times \rk(A)}$, $\U\!=\![U_1 \cdots U_N]$, and $\M\!=\!\tfrac{1}{N} \text{blkdiag}(M, \dots, M)$.
By Woodbury:
$(\nabla^2_p f_N)^{-1}\!=\!\tfrac{1}{\lambda} I \!-\! \tfrac{1}{\lambda^2} \U ( \M^{-1} \!+\! \tfrac{1}{\lambda} \U^\top \U )^{-1} \U^\top$.
This inverts a system of size $N\rk(A)$.
Since $N$ scales as $\cO(\ln|\cR|)$ (Prop.~\ref{prop:concentration_combined}) and $\rk(A) \!\propto\! |\cE|$, when $N|\cE| \!\ll\! |\cR|$ this reduces complexity from $\cO(|\cR|^3)$ to $\cO((|\cE| \ln |\cR|)^3)$.

\subsection{Details of Exp4 on GEANT}\label{sec:app:GEANT}
$K$ is generated using the shortest $k_{\text{paths}}\!=\!3$ candidate routes per OD pair,
$\mu \!=\!0.85x_{\max}$, $\sigma\!=\!1$, $\rho\!=\!0.05$,
$N\!=\!5000$, $p_{\min}\!=\!\zeros$, $p_{\max}\!=\!100 \cdot \ones$,
$\lambda\!=\!10^{-2}$.
Initial price $p_0\!=\!B^\dagger(\tfrac{x_{\max}}{2}-\mu)$.
Demand 
$l\!=\!\alpha K \Avg[\Pi_{[\zeros,x_{\max}]}Bp_0+\zeta]$ with $\alpha \!\sim\! \cU(0.3,0.7)$.
We report the DZ fraction at a point $p$ as
$\IE\big[\II(Bp+\zeta \notin\! [\zeros, x_{\max}])\big]$.
For TRSQP, we set 300 maximum iterations with tolerance 
$10^{-6}$ and constraint tolerence $10^{-3}$.
$\Delta_0=10, \Delta_{\max}=50, \Delta_{\min}=10^{-6}, \nu_0=1$.

\subsection{Details of Exp5 on Twitch social network}\label{sec:app:subsec:exp5}
We use the Stanford SNAP dataset on Twitch ENGB Networks\footnote{\url{https://snap.stanford.edu/data/twitch-social-networks.html}} . The data generation is similar to Exp4: Diagonal entries $B_{rr}\sim\cU [-3,\,-1.5]$;
cross-elasticities are nonzero only within the same OD pair,
drawn from a scaled $\mathrm{Beta}(2,8)$ distribution on $[0.3,\,0.8]$
and rescaled per row to ensure strict diagonal dominance
($|B_{rr}| \!>\! \textstyle\sum_{r'\neq r}|B_{rr'}|$).
$K$ is generated using the shortest $k_{\text{paths}}=3$ candidate routes per OD pair, $\mu \in \texttt{rand} \cdot x_{\max}$ with $\texttt{rand}\in [0.7, 1]$, $\sigma\!=\!1.5$, $\rho\!=\!0.05$,
$N\!=\!5000$, $p_{\min}\!=\!\zeros$, $p_{\max}\!=\!100 \cdot \ones$,
$\lambda\!=\!10^{-2}$. With initial price $p_0\!=\!B^\dagger(\tfrac{x_{\max}}{2}-\mu)$, the demand 
$l\!=\!\alpha K \Avg[\Pi_{[\zeros,x_{\max}]}Bp_0\!+\!\zeta]$ with $\alpha \!\sim\! \cU(0.3,0.7)$. 
For TR-SQP: 10 iterations, tolerance $10^{-6}$, constraint tolerence $10^{-3}$.
$\Delta_0\!=\!10, \Delta_{\max}\!=\!50, \Delta_{\min}\!=\!10^{-6}, \nu_0\!=\!1, \rho_{\text{th}}\!=\!10^{-2}$,
$\mu_V\!=\!
\max\!\big(\lambda\|p_{\max}\|_2\|\|B^{-1}\|_\infty, 10^{-6}\big)$.

Experiments ran on AMD EPYC 7413: 16 cores, 141 GB RAM, swarm a100 partition. We use 4 parallel workers for the sparse solver.

\subsection{Additional Experiments 4b}\label{sec:exp:subsec:nguyen}
We validate ONP on the Nguyen‑Dupuis synthetic EV charging benchmark%
\footnote{\url{https://github.com/bstabler/TransportationNetworks/tree/master/Nguyen-Dupuis}}, a widely used small-scale traffic assignment testbed modelling a urban corridor with $|\cV|=21$, $|\cE|=56$ directed edges, $|\cK|=4$ and $|\cR|=12$.
All other experimental settings ($B$ generation, initialization, parameters) follow  \cref{sec:exp:subsec:geant}.
 
\begin{figure}[ht!]
\begin{minipage}{0.32\linewidth}
    $\vcenter{\hbox{\includegraphics[width=\linewidth]{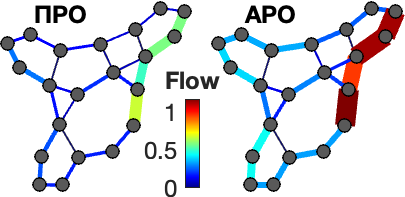}}}$
\end{minipage}
\begin{minipage}{0.59\linewidth}
    $\vcenter{\hbox{%
        \small
        \begin{tabular}{l c c}
            Metric & $\Pi$PO (ours) & APO \\
            \hline
            Final cost & $8.11$ & $16.36$ \\
            max constraint violation & $1.00\!\times\!10^{-2}$ & $3.73\!\times\!10^{-3}$ \\
            DZ fraction at solution & $89.1\%$ & $85.7\%$ \\
            Mean $\bar\pi_r$/Min $\pi_{\min}$ & $0.109/0.046$ & --- \\
            $\mathrm{nnz}(A)/(|\cR||\cE|)$ & $9.52\%$ & --- 
        \end{tabular}%
    }}$
\end{minipage}
\caption{Results on the Nguyen‑Dupuis EV network (mean path length $5.33$).
Robust initialisation sets the mean expected flow to $48.75\%$ of capacity with initial DZ fraction $85.8\%$.}
\label{fig:nguyen}
\end{figure}

\noindent
$\Pi$PO reduces the cost by $50.4\%$ compared to APO.
The APO solution is far from stationarity (the gradient mismatch $\|\nabla f_{\Pi\mathrm{PO}}(p_{\mathrm{APO}})\| = 68.46$), confirming that omitting the Jacobian severely biases the optimum even on a smaller network.

\noindent
At the $\Pi$PO solution $89.1\%$ of routes are dead‑zoned (APO: $85.7\%$).
The algorithm encountered DZ escape events.
The minimum activation probability $\pi_{\min}=0.046$.
All these shows that the escape mechanism is essential.

\paragraph{Sioux Falls}
We run the same experiment on the Sioux Falls ($|\cR|=1528$, see Fig.\,\ref{fig:SiouxFalls}), a canonical mid-scale traffic benchmark representing the road network.
$\Pi$PO improves the cost over APO for 74.35 $\pm$1.36\%.
\begin{figure}[h!]
\centering
\includegraphics[width=0.6\linewidth]{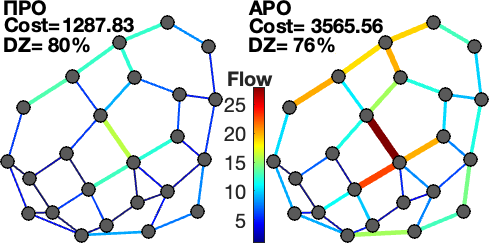}
\caption{A typical result on the Sioux Fall network ($|\cV|=24, |\cE|= 76$, $|\cK|=538$).}
\label{fig:SiouxFalls}
\end{figure}
\end{document}